\documentstyle[12pt]{article}

\newtheorem{defn}{Definition}[section]
\newtheorem{lem}[defn]{Lemma}
\newtheorem{thm}[defn]{Theorem}
\newtheorem{prop}[defn]{Proposition}
\newtheorem{cor}[defn]{Corollary}
\newtheorem{rem}[defn]{Remark}
\newtheorem{ex}[defn]{Example}
\newtheorem{assu-nota}[defn]{Notation - Assumption}
\newtheorem{assu}[defn]{Assumption}

\newtheorem{nota}[defn]{Notation}

\newcommand{\C}{{\bf C}}
\newcommand{\epsi}{\epsilon}
\newcommand{\pp}{{\bf P}}

\newcommand{\proof}{{\bf Proof:}\,\,}
\newcommand{\inv}{^{-1}}

\newcommand{\qed}{$\quad \diamond$\par\smallskip}
\newcommand{\OO}{{\cal O}}

\newcommand{\OZ}{{\cal O}_Z}

\newcommand{\Si}{\Sigma}

\newcommand{\Proj}{{\rm Proj}}
\newcommand{\Z}{{\bf Z}}
\newcommand{\spec}{\mbox{spec}}

\newcounter{rom}
\newcounter{arab}
\newcounter{alpha}

\title{Prym varieties and the canonical map of surfaces of general type}
\author{by Ciro Ciliberto, Rita Pardini and Francesca Tovena }
\date{}

\begin{document}
\maketitle
\begin{abstract} 
Let $X$ be a smooth complex surface of general type and let $\phi:X\to
{\bf P}^{p_g(X)-1}$ be the canonical map of $X$. Suppose
that the image $\Sigma$ of $\phi$ is a surface and that $\phi$ has
degree $\delta\geq 2$. Let $\epsilon\colon S\to \Sigma$ be a desingularization
of $\Sigma$ and assume that the geometric genus of $S$ is not zero.  Beauville
(\cite{Beau}) proved   that in this case the  surface $S$ is of general type and
$\epsilon$ is the canonical map of $S$. Beauville also constructed the only
infinite series of 
examples $\phi:X\to \Sigma$ with the above  properties 
that was up to now available in the literature. This construction has lead us to
introduce the notion of a {\em good generating pair}, namely of a pair $(h:V\to
W, L)$ where $h$ is  a  finite morphism of surfaces and $L$ is a nef and big line
bundle of $W$ satisfying  certain assumptions. The most important of these  are:
i) $|K_V+h^*L|=h^*|K_W+L|$, and ii) the general curve $C$ of
$L$ is smooth and non-hyperelliptic. We show that, by means of a construction
analogous to the
one of Beauville's,  every good generating pair gives rise to an infinite series
of surfaces of
general type whose  canonical map is $2$-to-$1$ onto a canonically embedded
surface. In this way we
are able to construct more infinite series of such surfaces (cf. section $3$). 
In addition, we
give bounds on the invariants of good generating  pairs and show that there exist
essentially only $2$ good  generating pairs with
$\dim |L|>1$. The key fact that we exploit for obtaining these results is that
the Albanese variety $P$ of
$V$ is a Prym variety and that the fibre of the Prym map over $P$ has positive
dimension.
\end{abstract}

\section{Introduction}

Let $X$ be a smooth surface of general type and let
$\phi\colon X\to
\Sigma\subseteq {\bf P}^{p_g(X)-1}$ be the canonical map of $X$,  where
$\Sigma$ is the image of $\phi$. Suppose that
$\Sigma$ is a surface and that $\phi$ has degree $\delta\ge 2$. Let 
$\epsilon\colon S\to \Sigma$ be a  desingularization  of 
$\Sigma$. A classical result, which goes back to Babbage  \cite{bab}, and has
been more recently proved by Beauville, \cite{Beau} (see also \cite{babbage}),
says that either $p_g(S)=0$ or $S$ is of general type and $\epsilon\colon S\to
\Sigma$ is the canonical map of $S$. In the latter case we have a 
dominant rational map $\psi\colon X\to S$ of degree $\delta$, which we call a
{\it good canonical cover} of degree $\delta$ (see definition \ref{ct} for a
slightly more general definition).\par

While there is no problem at all in exhibiting as many examples as one likes of
the former type, i.e. where $p_g(S)=0$ (see \cite{Beau}), not so many good
canonical covers are
available in the current literature.  In few sporadic examples of such covers
 the surface $X$ is {\it regular}
 (see \cite{zagier}, \cite{Beau}
proposition 3.6, \cite{babbage} theorem 3.5, \cite{enriques}, \cite{supcan}). On
the other hand, there is an
interesting construction, due to Beauville (see \cite{Ca}, 2.9 and \cite{MP}),
which produces an infinite series of such covers of degree $2$ where $X$ is {\it
irregular}, precisely of irregularity $2$. Beauville's construction is recalled
in \S 4 and in example \ref{Beauville}. The resulting canonical covers have been
extensively studied in \cite{MP},
where they have been classified in terms of their birational invariants. \par

In our attempts to find more examples of canonical covers, we have been lead to 
understand Beauville's construction better. In particular we extracted from it
its main
features, and this lead us to give a definition, the one of a {\it good
generating pair}
(see \ref{bpair} for a more general definition), which, roughly speaking, is the
following. A good generating pair $(h\colon V\to W, L)$ is the datum of a finite
morphism $h\colon V\to W$ of degree $2$ between surfaces, $V$ smooth and irreducible, $W$
with isolated double points of type $A_1$, and $L$ a nef and big line bundle on $W$.
Furthermore one requires that $|L|$ has at least dimension $1$ and contains a smooth,
irreducible, non--hyperelliptic curve $C$, that $h^*K_W=K_V$ (this means that $h$ has
only isolated
ramification points, corresponding to the double points of $W$) and that the pull-back
of the adjoint linear system $|K_W+L|$ is the complete linear system $|K_V+h^*L|$.\par

However cumbersome and un--motivated this definition may appear at a first glance, it
turns out to be rather useful for constructing canonical covers. Indeed one
finds many of these in the following way (see \S 3 for details).  Consider the map
$\tilde{h}=h\times  Id\colon V\times \pp^1\to W\times \pp^1$ and the projections $p_i$,
$i=1,2$,  of $W\times {\bf P}^1$ onto the two factors. A general surface
$\Sigma\in |p_1^*L\otimes p_2^*\OO_{\pp^1}(n)|$, $n\ge 3$, has only points of type $A_1$
as singularities. We set $X=\tilde{h}^*(\Sigma)$, $\phi=\tilde{h}|_X\colon  X\to
\Sigma$, $\epsilon\colon S\to \Sigma$ the minimal desingularization,
$\psi=\epsilon^{-1}\circ \phi$.
Then, using adjunction  both on $V\times \pp^1$ and $W\times\pp^1$, one sees that
$\psi:X\to S$ is a good canonical cover of degree $2$.\par

General properties of generating pairs $(h\colon  V\to W, L)$ are studied in \S
\ref{generalita} (see also \S \ref{altricasi}, where some information about {\it higher
degree} generating pairs has been collected). In particular, we see that $V$ and $W$
have the
same Kodaira dimension (see proposition \ref{kappa}), and, while $W$ is always regular,
$V$,instead, is irregular, and its irregularity can be expressed in terms  of the genus
$g$ of the general curve
$C\in |L|$ and of the degree of $h$ (see proposition \ref {regular}). We also give
formulas for the invariants of the canonical covers arising from a given generating pair
(see proposition \ref{explain}). \par

It is interesting to notice that Beauville's example is essentially characterized by the
fact that $V$ and $W$ have Kodaira dimension $\kappa=0$ (see proposition \ref{kappa0}
for a
more precise statement). The case of Kodaira dimension $1$ is also rather restricted, as
proposition \ref{kappa1} shows.\par

Beauville's example corresponds to the case in which $V$ is a principally polarized
abelian surface, $W$ is its Kummer surface, and $L$ is the polarization on $W$ which
lifts to a symmetric principal polarization on $V$. Unfortunately, more generating pairs
do not easily show up. The only ones which we know about are listed in section \S
\ref{esempi}. These  give rise to more infinite series of good canonical
covers which wait for a deeper understanding, like, as we said, in \cite{MP} has been
done for Beauville's examples.\par

The difficulty in finding generating pairs is not casual. This is explained in \S
\ref{grado2}, and this is where {\it Prym varieties} come into the picture. If $(h\colon 
V\to W, L)$ is a generating pair, and $C$ is a general curve in $|L|$, of genus $g$, then
$h^*C=C'\to C$ is an unramified double cover, with a related Prym variety $Prym(C',C)$.
In theorem \ref{prym} we prove that $Prym(C',C)$ is naturally isomorphic to the Albanese
variety of $V$. As a consequence we find that, if the generating pair is good,
the Albanese image of $V$ is a surface and therefore the Kodaira dimension of
$V$ and $W$ is non--negative (see corollary \ref{albsurface}). Moreover, some general
facts about irregular surfaces and isotrivial systems of curves on them, which have been
collected in \S \ref{irrfibr}, imply that the Prym map has an infinite fibre at the
cover $h\colon  C'\to C$ (see proposition
\ref{notfinite}). This, together with results about the fibre of the Prym map
due to several authors (see \S \ref{grado2} for references), enable us to
prove that, if $(h\colon V\to W, L)$ is a good pair, then  one has the bounds
$g\leq 12$ for the genus $g$ of $C$ and
is $q\leq 11$ for the irregularity $q$ of $V$ (see theorem \ref{g<12} and proposition
\ref {genere}). We suspect that, along the same lines, it should be possible to improve
this bound for $g$ and $q$, but this would preliminarly require a deepening of our
understanding of the fibres of the Prym map. For instance we would like to
know answers to questions like: when may these fibres contain rational curves?
Problems, of course, of independent interest.\par

Finally, using  Reider's method, we obtain the bound $L^2\leq 4$ (see proposition \ref
{Reider}), so that one really sees why there are not so many possibilities for a good
generating pair. We give  a complete classification of good pairs with
$L$ ample and $h^0(L, W)>2$.\par 

These satisfy $h^0(W,L)\leq 4$ and $L^2=3, 4$. The only
example with $h^0(W,L)=4$ is Beauville's one (see corollary \ref {classifica}).
The cases $h^0(W,L)=3$ and $L^2=3$ or
$L^2=4$ but $|L|$ with a base point are studied in \S \ref{sezgood}.
(see theorem \ref{trigonalcase}); we find that the former case corresponds
either to example \ref {Rita} or  to a suitable modification of Beauville's
example, while the latter does not occur. These cases share the feature that
the general curve $C$ is trigonal, and we take advantage, in the proof of our
classification theorem \ref{trigonalcase}, of a globalization to $V$ of the
well known {\it trigonal construction} (\cite{cas2}), which is the inverse of the
equally famous {\it Recillas' construction}. A different proof of the same
result is sketched in  remark \ref {trix}.\par

Finally we prove that $L^2=4$, $h^0(W,L)=3$ does not occur (see corollary
\ref {quattrocor}).   In the {\it pencil case} $h^0(W,L)=2$ (in which  there
are  examples, like \ref{Ciro}, but a classification is still lacking) we
show that the possibility $L^2=4$ is severely
restricted  (see
corollary \ref {quattrocor}).\par 

Using similar ideas,
we are able to construct an infinite family of good canonical covers with $X$ regular. We
will be back on this in a forthcoming paper.

\noindent{\bf Notation and conventions:} all varieties are defined over
the field of complex numbers. A map between varieties is a rational map, while a
morphism is a rational map that is regular at every point. We do not
distinguish between Cartier divisors and line bundles and use the additive and
multiplicative notation interchangeably. The Kodaira dimension of a variety
$X$ is denoted by $\kappa(X)$. We denote by $\sim_{num}$ the numerical
equivalence between divisors on a smooth surface.

\section{Canonical covers and generating pairs}\label{definizioni} 

\begin{nota}\label{pgq} Let $S$ be a surface with canonical singularities, i.e.
either smooth or with rational double points, so that in particular $S$ is Gorenstein.
We denote by $K_S$ the canonical divisor of $S$, and we let
$p_g(S)=h^0(S,K_S)=h^2(S,\OO_S)$ be the  {\em geometric genus} and $q(S)=h^1(S, \OO_S)$
the {\em irregularity}. If $p_g(S)\geq 2$, the {\em canonical map} of $S$ is the
rational map
$\phi\colon  S\to \pp^{p_g(S)-1}$ defined by the moving  part of the
canonical system $|K_S|$ of $S$. If $S_0$ is the open set of smooth points of
$S$ and  $\epsi\colon S'\to S$ is any desingularization, then $p_g(S)=p_g(S')$
and  $q(S)=q(S')=h^0(S', \Omega^1_{S'}) = h^0(S_0, \Omega^1_{S_0})$. The
Albanese map  of $S'$ factors through $\epsi$, since the exceptional locus of
$\epsi$ is a union of rational curves, and so we can speak of the Albanese map
of $S$. \end{nota}

Let $X$ be a smooth surface of general type and let $\phi\colon  X\to
\Sigma\subseteq {\bf P}^{p_g(X)-1}$ be the canonical map of $X$, where
$\Sigma$ is the image of $\phi$. We assume that $\Sigma$ is a surface and that
$\phi$ has degree $d\ge 2$, and we denote by  $\epsilon\colon  S\to \Sigma$
a desingularization of $\Sigma$. We recall the following theorem due to
Beauville, \cite{Beau}, Thm. 3.4.

\begin{thm} \label{beau} Under the above assumptions,  either:\par
\noindent (i) $p_g(S)=0$ or;\par
\noindent (ii) $S$ is of general type and $\epsilon\colon  S\to
\Sigma$ is the canonical map of $S$.
\end{thm}

\noindent
We introduce  some terminology for surfaces verifying condition (ii) of Theorem
\ref{beau}:

\begin{defn}\label{ct} Let $X, S$ be smooth surfaces of general type. Let
$\psi\colon  X\to S$ be a dominant rational map of degree $d\ge 2$. Assume
that: \begin{list} {(CC\arabic{arab})}{\usecounter{arab}}
\item $p_g(X)=p_g(S)$;
\item the canonical image  of $S$ is a surface $\Si$.
\end{list}
In this case the canonical map $\phi:X\to \Si$ of $X$ is the composition of
$\psi$ and the canonical map $\epsilon\colon  S\to \Sigma$ of $S$, and we say
that $\psi\colon  X\to S$ is a  {\em canonical cover of degree $d$}. If $\epsilon\colon
S\to \Sigma$ is birational, then we say that the canonical cover is {\em good}.\par
 \end{defn}
A few sporadic examples of canonical covers are available in the literature
(\cite{zagier},\cite{Beau} prop. 3.6, \cite{babbage} thm. 3.5, \cite {enriques},
\cite{supcan}). However, so far, there is only one construction, due to
Beauville (see \cite{Ca}, 2.9 and \cite{MP}), which produces an infinite series of such
covers. We  recall it next.\par 
Let $V$ be a principally polarized abelian surface  such that the principal
polarization $D$ is irreducible, and let $h\colon V\to W$ be the quotient map
onto the Kummer surface $W=V/<-1>$. The surface $W$ can be embedded into
$\pp^3$ as a quartic surface via a complete linear system $|L|$ such that
$h^*|L|=|2D|$. Consider the map $\tilde{h}=h\times Id\colon V\times \pp^1\to
W\times \pp^1$ and the projections $p_i$, $i=1,2$, of $W\times {\bf P}^1$
onto the two factors. A general surface $\Sigma\in |p_1^*L\otimes
p_2^*\OO_{\pp^1}(n)|$, $n\ge 3$, has only points of type $A_1$ as
singularities. We set $X=\tilde{h}^*(\Sigma)$, $\phi=\tilde{h}|_X\colon  X\to
\Sigma$, $\epsilon\colon  S\to \Sigma$  the minimal desingularization,
$\psi=\epsilon^{-1}\circ \phi\colon  X\to S$. Then it is easy to check, using
adjunction both on $V\times \pp^1$ and $W\times\pp^1$, that $ \psi\colon  X\to
S$ is a good canonical cover of degree $2$.\par
We wish to study to what extent this construction can be generalized. We
introduce  a class of pairs $(h\colon V\to W,L)$, where $h\colon V\to W$ is a
finite morphism of surfaces and $L$ is a line bundle on $W$, in such a way
that by applying the above construction to $(h\colon V\to W,L)$ one gets an
infinite series of canonical covers.\par
\begin{defn}\label{bpair} Consider a pair $(h\colon V\to W, L)$, where $h$ is
a finite morphism   of degree $d\geq 2$ between irreducible surfaces, $V$
smooth, $W$ with at most canonical singularities
and $L$ is a line bundle on $W$, such that:
\begin{list} {(GP\arabic{arab})}{\usecounter{arab}}
\item $K_V=h^* K_W$;
\item $h^0(W,L)\geq 2$ and  $L$ is big, i.e. $ L^2>0$;
\item  the general
curve $C$ of   $|L|$ is smooth of genus $g\ge 2$ and the curve $C':=h^*C$ is not
hyperelliptic;
\item $p_g(V)=p_g(W)$, $h^0(V,K_V+h^*L)=h^0(W,K_W+L)>0$.
\end{list}
We call $(h\colon V\to W, L)$  a {\em degree $d$ and genus $g$ generating
pair} of canonical covers, and we denote by $L'$ the line bundle $h^*L$ on
$V$.
The pair is said to be {\em minimal} if both $V$ and $W$ are minimal.\par
The generating  pair is called {\em good} if the general $C$ of
$|L|$ is not hyperelliptic (hence $g\ge 3$ in this case). \par
\end{defn}
Notice that condition (GP1) is equivalent  to the fact that $h$ is ramified
only over the singular points of $W$. Condition (GP3) and Bertini's theorem
imply that the general curve $C$ in $|L|$ is smooth and irreducible, hence $L$
is also nef, i.e. $L D\geq 0$ for every effective
divisor $D$ on $V$. The assumption that $C'$ is not hyperelliptic is a technical
condition whose meaning will be clearer later (cf. for instance theorem \ref{prym}).
Finally, the base points of $|L|$, if any, are smooth points of $W$. \par In the rest of
this section we show that by applying the original construction of Beau\-ville, to a
(good) generating pair one obtains an infinite series of (good) canonical covers, and we
compute the invariants of such canonical covers. In order to do this, we need the
following result, that will be proven later (cf. proposition \ref{regular}):
\begin{prop}\label{regular1} If  $(h\colon V\to W, L)$ is  a generating pair,
then $q(W)=0$. \end{prop} 
\medskip We introduce now some more notation:
\begin{nota}\label{notbpair} Given a generating pair  
$(h\colon V\to W, L)$ of degree $d$ and genus $g$, we
denote by $p_i, i=1,2$, the projections of $W\times \pp^1$ onto the two
factors and we write $\tilde{h}=h\times Id\colon V\times\pp^1\to 
W\times\pp^1$. We denote by ${\cal L}(n)$ the line bundle
$p_1^*L\otimes p_2^*\OO_{\pp^1}(n)$, where $n$ is a positive integer.\par 
In addition, we let $\Sigma\in |{\cal L}(n)|$ be a general surface,
$Y=\tilde{h}^*(\Sigma)$. We denote by $\epsilon\colon  S\to \Sigma$ and 
$\epsilon'\colon X\to Y$ the minimal desingularizations,  by
$f$ the map
$\tilde{h}|_X\colon  X\to \Sigma$, and by 
$\psi$ the map $\epsilon^{-1}\circ f\circ \epsilon'\colon  X\to S$. \end{nota}

\begin{prop}\label{explain}  We use notation \ref{notbpair}.\par \noindent Let $(h\colon
V\to W, L)$ be a generating pair of degree
$d$ and genus $g$. If $n\geq 3$, then one has:\par
i) let   $\Sigma\in |{\cal L}(n)|$ be general and let $Y=h^*\Si$; $\Si$ and $Y$ are
surfaces  of general type with at most canonical singularities.  If in addition $L$
is ample, then   $S$ and
$X$  are both minimal;\par ii) $\psi\colon  X\to S$ is a   canonical cover of degree
$d$, that is said to be {\em $n$--related} to the generating pair $(h\colon V\to W,
L)$. If the generating pair is  good, then $\psi\colon  X\to S$ is a good canonical
cover, while if the generating pair is not good then  the canonical map of
$S$ is
$2$-to-$1$ onto  a rational surface;\par
iii) the  invariants of $S$ are: $p_g(S)= np_g(W)+ (n-1)g$, $q(S)=0$,
$K_S^2=n(K_W^2-L^2)+8(n-1)(g-1)$;\par
iv) the invariants of $X$ are: $p_g(X)=p_g(S)$ and 
$q(X)=(d-1)(g-1)$, $K^2_X=d\,K_S^2$.\par
\end{prop} \proof  Recall that by condition (GP2) of definition
\ref{bpair} the general curve of $|L|$ is smooth, and thus, in particular,
$|L|$ has no fixed components. Thus also the linear system
$|{\cal L}(n)|$ has no fixed components and is not composed with a pencil.
Therefore its general member $\Sigma$ is irreducible. Moreover the set of base
points of $|{\cal L}(n)|$ is the inverse image via $p_1$   of the set of base
points of  $|L|$ and thus it is a finite union of fibres of $p_1$. Using
Bertini's theorem and the fact that the general curve of
$|L|$ is smooth, one proves that the singularities of the general 
$\Sigma\in|{\cal L}(n)|$ at points  of the fixed locus of $|{\cal L}(n)|$ are
finitely many rational double points of type $A_r$. Now, the projection $p_1$
restricts to a generically finite map $p\colon \Sigma\to W$ of degree $n$ 
which,  by Bertini's theorem again, is unramified over the singular points of
$W$. So the  general $\Sigma$ has, over each singular  point $x$ of $W$, $n$
singularities which are analytically equivalent to the one $W$ has in 
$x$ (i.e. $n$ canonical singularities) and it is smooth at points that are
smooth for $W\times\pp^1$ and are not base points of $|{\cal L}(n)|$. To
describe the singularities of
$Y=\tilde{h}^*(\Si)$, we notice that the restriction $Y\to\Sigma$ of
$\tilde{h}$ is ramified precisely over the singularities of $\Si$ that occur at 
singular points of $W\times\pp^1$; so $Y$ has $d$ singularities analytically
isomorphic to those of $\Si$ over  each of those singular points of $\Si$ that
occur at base points of $|{\cal L}(n)|$ and it is smooth elsewhere, since it
is general in $\tilde{h}^*|{\cal L}(n)|$. In conclusion the singularities of
$Y$ and $\Si$ are canonical, and their invariants, which we  now compute,  
are equal to those of $X$,  $S$, respectively.\par
By the adjunction formula and condition (GP1) in definition \ref{bpair}, one
has $K_\Sigma=(K_{W\times\pp^1}+\Sigma)|_\Sigma= (p_1^*K_W+
{\cal L}(n-2))|_\Sigma$ and
$K_X=(K_{V\times\pp^1}+X)|_X=\tilde{h}^*(K_{W\times\pp^1}+\Sigma)|_X=
\psi^*(K_ \Sigma)$, and thus $K_S^2=K_\Sigma^2=n(K_W^2-L^2)+8(n-1)(g-1)$,
$K_X^2=dK_S^2$ . To compute the remaining invariants of $S$, $\Sigma$ and $X$,
one considers the long cohomology sequences associated to the restriction
sequences  $$0\to K_{W\times\pp^1}\to K_{W\times\pp^1}+{\cal L}(n)\to 
K_\Sigma\to 0$$ and $$0\to K_{V\times\pp^1}\to
K_{V\times\pp^1}+\tilde{h}^*{\cal L}(n)\to  K_X\to 0.$$ 
By Kawamata-Viehweg's vanishing theorem, we have $h^i(W\times \pp^1,
K_{W\times\pp^1}+{\cal L}(n))=h^i(V\times \pp^1,
K_{V\times\pp^1}+\tilde{h}^*{\cal L}(n))=0$ for $i>0$. Hence:
$$p_g(S)=h^0(\Sigma,K_\Sigma)=$$
$$=h^0(W\times\pp^1, K_{W\times\pp^1}+{\cal L}(n))+h^1(W\times\pp^1, 
K_{W\times\pp^1})=$$
$$=h^0(W, K_W+L)(n-1)+p_g(W)=np_g(W)+(n-1)g$$
where the last equality follows again from Kawamata--Viehweg's vanishing and
the last equality but one follows from $q(W)=0$. Therefore, by the definition
of a generating pair: $$p_g(X)=h^0(V, K_V+L)(n-1)+p_g(V)=p_g(S).$$
\noindent A similar computation gives
$q(S)=q(W)=0$, $q(X)=q(V)=(d-1)(g-1)$.\par
The linear  system  $|K_\Sigma|$ contains the restriction of the system
$|p_1^*K_W+{\cal L}(n-2)|$, whose fixed locus is the inverse image via 
$p_1$ of the fixed locus of $|K_W+L|$. Let $C\in |L|$ be a smooth curve;
since $W$ is regular, the linear system $|K_W+L|$ restricts to the complete
canonical system $|K_C|$.  Thus $C$ does not contain any base point of
$|K_W+L|$. If $L$ is ample, this implies that $|K_W+L|$ has a finite
number of base points, none of which is also a base point of $|L|$. Thus in
this case the fixed locus of $|p_1^*K_W+{\cal L}(n-2)|$ intersects the general
$\Sigma$ in a finite number of points and, a fortiori, the canonical system of
$\Si$ has no fixed components and the surfaces $X$, $S$ are minimal.\par
Notice now that $|p_1^*K_W+{\cal L}(n-2)|$ separates the fibres of
$p_2|_{\Si}$,  since $n\ge3$. A  fibre $F$ of $p_2|_\Sigma$ is identified by
$p_1$ with a curve $C\in |L|$ and the restriction of $|p_1^*K_W+{\cal L}(n-2)|$
to $F$ is identified with the restriction of $|K_W+L|$ to $C$, which is the
complete canonical system $|K_C|$, since $W$ is regular. Thus, if the general
$C$ is not hyperelliptic, then the canonical map of $S$ is birational
and $\psi\colon  X\to S$ is a good canonical cover, while if the general
$C$ is  hyperelliptic then the canonical map of $S$ is of degree $2$ onto a
rational surface and  $\psi\colon  X\to S$ is a non good  canonical
cover. \hfill\qed

\medskip
Since we aim at a classification of generating pairs, we find useful to
introduce a notion of blow-up. We will show (cf. corollary \ref{mintame}) that in
most  cases that almost every generating pair is obtained from a minimal one by a
sequence of blow-ups.
\begin{defn}\label{birational} Let $(h\colon  V\to W, L)$ be a generating
pair of degree $d$ and genus $g$. Let $x\in W$ be a smooth point. Then we can
consider the cartesian square:
$$\begin{array}{rcccl}
\ &V'&\rightarrow &V&\ \\
\scriptstyle{h'}\!\!\!\!\!\! & \downarrow &  &\downarrow &\!\!\scriptstyle{h}\\
\ &W' & \stackrel{f}{\rightarrow} & W &
\end{array}$$
\noindent where $f\colon  W'\to W$ is the blow-up of $W$ at $x$, with
exceptional divisor $E$ and, accordingly, $V'$ is the blow-up of $V$ at the
$d$ points $x_1,...,x_d$ of the fibre of $h$ over $x$.  Fix $m=0$ or $1$ and assume
that:\par
\noindent (i) $L^2>m^2$;\par
\noindent (ii) $h^0(W', f^*L-mE)\geq 2$ and the general curve $C\in |
h'^*L-mE|$ is smooth. \par
Then the pair $(h'\colon  V'\to W',f^*L-mE)$ is again a generating pair. We  say that
it is obtained from $(h\colon  V\to W, L)$ by a simple blow-up. 
The blow-up is said to be {\em essential} if $m=1$ and {\em inessential} if $m=0$.
\end{defn}
The reason why we only consider $m\le 1$ in the above definition is that generating
pairs satisfy the inequality $L^2\le 4$ (cf. prop. \ref{Reider} and prop. \ref{d=3}).

\section{Examples of generating pairs}\label{esempi}

In this section we describe some examples of generating pairs. 

\begin{ex} \label{Beauville} Beauville's example. {\rm (see 
\cite{Ca}, 2.9,
\cite{MP}, example 4 in section 3). This example has already been described in
section \ref{definizioni}: $V$ is a principally polarized abelian surface with an
irreducible polarization $D$, $W$ is the Kummer surface of $V$, $h\colon  V\to W$ is the
projection onto the quotient, and $L$ is an ample line bundle on
$W$ such that the class of $L'=h^*L$ is equal to $2D$.  This generating pair is
good and
therefore so is any related canonical cover. More precisely, by proposition 
\ref{explain}, an $n$-related canonical cover $\psi\colon  X\to S$ is minimal,
with geometric genus $4n-3$. The invariants of $S$ and $X$ satisfy the relations:
$$\quad K^2_S=3p_g(S)-7;\quad K^2_X=6p_g(X)-14; \quad q(X)=2.$$
According to \cite{MP}, Thm. 4.1, this is the only good generating  pair  such
that the related canonical covers satisfy 
$K^2_X=6p_g(X)-14$ and $K_X$ is ample.\par

Notice that, if, in the above situation, the polarization $D$ on $V$ is not irreducible,
then the same construction produces a generating pair which is no
longer good (cf. also \cite{MP}, example 2 in section 3). We will
refer to this example as to the {\it non good Beauville's
example}.} \end{ex}

\begin{ex} \label{Ciro} A good generating pair of degree $2$ and genus $3$.
{\rm (cf. also \cite{bic},  example (c), page 70). Let $A$ be an abelian
surface with an irreducible  principal polarization $D$,  let $p\colon V\to A$
be the double  cover branched on a symmetric divisor $B\in |2D|$ and
such that $p_*\OO_V=\OO_A\oplus \OO_A(-D)$. Since $K_V=p^*(D)$, the invariants
of the smooth surface $V$ are: $p_g(V)=2$, $q(V)=2$, $K^2_V=4$. By the symmetry
of $B$,  multiplication by $-1$ on $A$ can be lifted  to an involution $i$ of $V$ that
acts as the identity on $h^0(V, K_V)$. We denote by $h\colon V\to W=V/<i>$ the
projection onto the quotient. We observe that
$p_g(W)=2$, $q(W)=0$, $K^2_W=2$ and the only singularities  of the surface $W$ are $20$
ordinary double points. In addition,
$h^0(W,2K_W)=\chi(\OO_W)+K^2_W=4=h^0(V, 2K_V)$, so that the bicanonical map of
$V$ factors through $h\colon V\to W$. An alternative description of $W$ is as
follows. One embeds, as usual, the Kummer surface $Kum(A)$ of $A$ as a quartic
surface in $\pp^3=\pp(H^0(A,2D)^\star)$. The surface $W$ is a double  cover of
$Kum(A)$  branched over the smooth plane section $H$ of $Kum(A)$ corresponding
to $B$ and on $6$ nodes (corresponding to the six points of order 2 of $A$
lying on $D$). The ramification divisor $R$ of $W\to Kum(A)$ is a canonical
curve isomorphic to $H$, and thus it is not hyperelliptic. This completes the
proof that $(h\colon V\to W, K_W)$ is a good generating pair. Notice that,
under suitable generality assumptions, $K_W$, as well as $K_V$, is ample. An
$n$-related canonical cover $\psi\colon  X\to S$ has geometric genus $5n-3$
and is, in general, minimal. The invariants of $S$ and $X$ satisfy the
relations: \begin{equation}\label{esempiociro}  5 K^2_S=16p_g(S)-32;\quad
5K^2_X=32 p_g(X)-64; \quad q(X)=2. \end{equation} }\end{ex}

\begin{ex}\label{Rita}  A good generating pair of degree $2$ and genus $4$ {\rm
(cf. \cite{bic},  example 3.13). Let $\Gamma$ be a non--hyperelliptic curve of genus
$3$ and let $V:=Sym^2(\Gamma)$. The surface $V$ is smooth minimal of general type with
invariants: $K^2_V=6$, $p_g(V)=q(V)=3$. If we embed  $\Gamma$ into $\pp^2$  via the
canonical system, then the canonical map of $V$ sends the  unordered pair $\{ p,q\}$ of
$V$ to the line $<p,q>\in\pp^{2*}$,  hence it is a degree $6$  morphism onto the plane.
There
is an involution $i$ on $V$ that maps $\{p,q\}\in V$ to $\{r, s\}$, where $<p,q>\cap
\Gamma=p+q+r+s$. The  fixed points of $i$ correspond to the $28$ bitangents of $\Gamma$
and the canonical map of $V$ clearly factors through the quotient map $h\colon V\to
W=V/<i>$.
Hence the  invariants of $W$ are: $p_g(W)=p_g(V)=3$,
$K^2_W= K^2_V/2=3$, $\chi(W)=(\chi(V)+7)/2=4$, and thus $q(W)=0$.  In addition we have
$h^0(W, 2K_W)=\chi(\OO_W)+K^2_W=7=h^0(V, 2K_V)$ and thus $|2K_V|=h^*|2K_W|$. In order to
complete the proof that $(h\colon V\to W, K_W)$ is a  good generating pair we remark 
that the general canonical curve $C$ of $W$ is  not hyperelliptic, since
the restriction of
$|K_W|$ to $C$ is a base--point free $g^1_3$. Notice that
$K_W$ and $K_V$ are ample. An $n$-related canonical cover $\psi\colon  X\to S$
is minimal,
of geometric genus $7n-4$. The invariants of $S$ and $X$ satisfy the relations:
\begin{equation}\label{esempiorita} 
7 K^2_S= 24 p_g(S)- 72 ;\quad 7 K^2_X= 48 p_g(X)- 144; \quad q(X)=3.
\end{equation}
An interesting question, concerning this example and the previous one, is
whether these are the only generating pairs such that the related canonical
covers have invariants
satisfying (\ref{esempiociro}) and (\ref{esempiorita}).} \end{ex}

\begin{ex} A non good generating pair of degree $3$ and genus $2$. {\rm Let
$p\colon W\to {\bf P}^2$ be the double cover of ${\bf P}^2$ ramified on an
irreducible sextic $B$ with 9 cusps ($B$ is the dual of a smooth cubic). The surface $W$
is a K3 surface whose singularities are 9 double points of type $A_2$. According to
\cite{BdF} (cf. also \cite{BiL}, \cite{barth}), there exists a smooth cover 
$h\colon V\to W$ of degree 3 ramified only at the 9 double points. The surface
$V$ is an abelian surface. Let $L=p^*({\cal O}_{{\bf P}^2}(1))$. Since $L$ is ample, we
have $h^0(W,L)=\chi(\OO_W)+\frac{1}{2}L^2=3=h^0(V, h^*L)$, and thus $(h\colon V\to W,
L)$ is a non good generating pair. An $n$-related, minimal, canonical cover $\psi\colon 
X\to S$ has geometric genus $4n-2$. The invariants of $S$ and $X$ satisfy the relations:
$$ K^2_S=2p_g(S)-4;\quad K^2_X=6p_g(X)-12;
\quad q(X)=2.$$ It is perhaps worth remarking that the surfaces $S$ thus
obtained have invariants lying on the Noether's line $K^2_S=2p_g(S)-4$. It
would be interesting to know whether there are
other canonical covers with so low geometric genus.} \end{ex} 

\begin{ex}\label{product}  A  series of non good generating pairs of degree $2$ with
unbounded invariants. {\rm For  $i=1,2$, let
$\phi_i\colon C_i\to\pp^1$ be a double cover, where $C_i$ is a smooth  curve of
genus $g_i>0$, and let $\sigma_i$ be the involution on $C_i$ induced by
$\phi_i$. We set $V=C_1\times C_2$, $W=V/<\sigma_1\times \sigma_2>$ and we
denote  by $h\colon V\to W$ the projection onto the quotient. We remark that
there exists a double cover $f\colon W\to \pp^1\times\pp^1$ such that
$\phi_1\times\phi_2\colon V\to\pp^1\times\pp^1$ factors as 
$\phi_1\times\phi_2=f\circ h$. We denote by $H$ a divisor of type $(1,1)$ on
$\pp^1\times\pp^1$ and we set $L=f^*H$.  Both systems $|K_V|$ and
$|K_V+h^*L|$ are clearly pull-back via
$\phi_1\times\phi_2\colon V\to\pp^1\times\pp^1$. This immediately implies that
$(h\colon V\to W, L)$ is a non good generating pair of degree $2$ and genus
$g_1+g_2+1$. One has: $p_g(W)=g_1g_2$, $q(V)=g_1+g_2$. An $n$-related canonical
cover $\psi\colon  X\to S$ has geometric genus $ng_1g_2+(n-1)(g_1+g_2+1)$ and
moreover $q(X)=g_1+g_2$.\par
Notice that, if $g_1=g_2=1$, we find again the non good Beauville's example (cf. example
\ref{Beauville}). } \end{ex}

\section{Auxiliary results on irregular surfaces}\label{irrfibr}
In this section we collect a few general facts on irregular surfaces that will be used
in the rest of the paper. We use notation \ref{pgq}.
\begin{prop}\label{albcurva} Let $h\colon V\to W$ be a finite morphism of
surfaces with
canonical singularities such that $K_V=h^*K_W$ and $p_g(V)=p_g(W)$. If 
$q(V)>q(W)>0$, then the Albanese image of $V$ is a curve. 
\end{prop} \proof The critical set $\Delta$ of $h$ is finite by assumption. We let
$W_0=W\setminus({\rm Sing} (W)\cup \Delta) $ and $V_0=h\inv  W_0$, so that the
restricted map $h\colon V_0\to W_0$ is a finite \'etale map between smooth surfaces. In
particular $h$ is flat, and there is a canonical vector bundle isomorphism
$h_*\OO_{V_0}\cong
\OO_{W_0}\oplus E$, where $h_*\OO_{V_0}$ and $E$ are locally free of ranks $d=\deg h$
and $d-1$ respectively.
Since $\Omega_{V_0}^i=h^*\Omega_{W_0}^i$, $i=1,2$,  one has
$h_*\Omega_{V_0}^i=\Omega_{W_0}^i\otimes h_*\OO_{V_0}=\Omega_{W_0}^i \oplus
(\Omega_{W_0}^i\otimes E)$. Notice that this decomposition as a direct sum is canonical. 
We set $M^i_+=H^0(W_0,\Omega^i_{W_0})$ and $M^i_-=H^0(W_0,\Omega_{W_0}^i\otimes E)$. We
deduce
that  $H^0(V_0,\Omega^i_{V_0})=H^0(W_0,h_*\Omega^i_{V_0}) =M^i_+ \oplus M^i_-$ and we
denote by $\pi^i_+$ the projection onto the first factor of this decomposition. To ensure
that the Albanese image is a curve, we show that $\tau_1\wedge \tau_2=0$ for
every choice
of $\tau_1, \tau_2 \in  H^0(V_0,\Omega^1_{V_0})=M^1_+\oplus M^1_-$. Noticing
that both
$M^1_+$ and $M^1_-$ are non-zero (since $q(V)>q(W)>0$), we only need to show that
$h^*\sigma \wedge \tau =0$ for every choice of $\sigma \in M^1_+$ and $\tau \in M^1_-$.
Indeed, to show that $\wedge^2 M^1_+=0$ we fix $(0\neq)\tau\in M^1_-$: if $\sigma_1,
\sigma_2\in M^1_+$, the vanishing $h^*\sigma_i\wedge \tau=0$ ($i=1,2$) means that 
$h^*\sigma_i$ is pointwise proportional to $\tau$ ($i=1,2$), so that $h^*\sigma_1$ and
$h^*\sigma_2$ are mutually pointwise proportional. Similarly one proves
that $\wedge^2 M^1_-=0$.\par

Since $p_g(V)=p_g(W)$,  $\pi^2_+$  is an isomorphism. Notice also that $\pi^i_+(h^*\sigma)=
\sigma$ for any $\sigma\in M^i_+$, and  that  $\pi^2_+(h^*\sigma\wedge\tau)=\sigma\wedge
\pi^1_+(\tau)$  for $\sigma\in M^1_+$ and $\tau\in  H^0(V_0,\Omega^1_{V_0})$.  Therefore
$h^*\sigma\wedge\tau=0$ for any $\sigma\in M^1_+$ and $\tau\in ker\, \pi^1_+=M^1_-$, as we
wanted. \hfill\qed\medskip

We recall the following results: 
\begin{prop}\label{isostruttura} (Serrano, \cite{Serrano}, section 1) Let $V$ be a smooth
surface, let
$C$ be a smooth curve,  and let $p\colon V\to C$ be an isotrivial fibration
with fibre $D$. Then there exist a curve $B$, a finite group $G$ acting both
on $B$ and  $D$, an isomorphism $f\colon C\to B/G$, and a  birational map
$r\colon V\to (D\times B)/G$, where $G$ acts diagonally on $D\times B$, such
that the following diagram commutes:

$$\begin{array}{rcccl}
\phantom{1} &V&\stackrel{r}{\rightarrow} &(D\times B)/G&
\phantom{1} \\
\scriptstyle{p}\!\!\!\!\!\! & \downarrow & \phantom{1} &
\downarrow &
\!\!\!\!\!\!\!\!\!\!\!\!\!\!\!\!\!\scriptstyle{p''}\\
\phantom{1} & C & \stackrel{f}{\rightarrow} & B/G &
\phantom{1}
\end{array}
$$  where $p''$ is the map induced by the projection $D\times B\to B$.
The irregularity $q(V)$ is equal to $g(C)+g(D/G)$. In particular, if
$q(V)>0$ and $g(C)=0$, then the Albanese image of $V$ is a curve isomorphic to
$D/G$ and the Albanese pencil is given by the composition $p'\circ r$, where
$p'$ is the map induced by the projection $D\times B\to D$.
\end{prop}

\begin{prop}\label{xiaofib} (Xiao, \cite{xiao}, Thm.1) Let $p\colon V\to \pp^1$ be a fibration
with fibres of genus $\gamma$. If
$p$ is not isotrivial, then $\gamma\ge 2q(V)-1$.\end{prop}

\medskip The next proposition combines the previous results.

\begin{prop}\label{xiaose} Let $V$ be a smooth surface with  a pencil 
$|D|$ such that the general curve $D$ of $|D|$ is smooth and
irreducible of genus $\gamma>1$; if the Albanese image of $V$ is a curve, then one (and 
only
one) of the following holds:\par

i) there exists a birational map $r\colon V\to D\times \pp^1$ such that
$D$ is the strict transform via $r$ of a fibre of the projection $D\times
\pp^1\to\pp^1$. In this case $\gamma=q(V)$;\par

ii) there exist an hyperelliptic curve $B$,  a free  involution $i$ on
$D$, and a birational map $r\colon V\to (D\times B)/\Z_2$, where $\Z_2$ acts on
$B$ as the hyperelliptic involution, on $D$ via $i$ and diagonally on $D\times
B$, such that $D$ is the strict transform via $r$ of a fibre of the projection
$(D\times B)/\Z_2\to B/\Z_2=\pp^1$. In this case $\gamma=2q(V)-1$.\par

iii) $\gamma>2q(V)-1$.\par

In particular, if $p$ is not isotrivial, then  iii) holds.
\end{prop} \proof Since the statement is essentially birational, up to blowing up the base
locus of $|D|$, we may assume that $D$ defines a morphism $p\colon V\to\pp^1$. Denote by
$\alpha\colon  V\to C$ the Albanese pencil. If $p$ is not isotrivial, then 
$\gamma\ge 2q(V)-1$ holds by proposition \ref{xiaofib}. If $\gamma=2q(V)-1$ ,
then by the Hurwitz formula the restriction of $\alpha$  to a smooth curve $D$ is an \'etale
cover of $C$, whose degree is $2$. Thus $p$ is isotrivial, contradicting the previous
assumption.\par

Assume now that $p$ is isotrivial. By  proposition \ref{isostruttura}, there is a
commutative diagram:
$$ \begin{array}{rcccl}
\phantom{1} &V&\stackrel{r}{\rightarrow} &(D\times B)/G&
\phantom{1} \\
\scriptstyle{p}\!\!\!\!\!\! & \downarrow & \phantom{1} &
\downarrow &
\!\!\!\!\!\!\!\!\!\!\!\!\!\!\!\!\!\scriptstyle{p''}\\
\phantom{1} & \pp^1 & \stackrel{f}{\rightarrow} & B/G &
\phantom{1}
\end{array}$$  
where $B$ is a curve, $G$ is a finite group acting on $B$ and on $D$ and acting
diagonally on $D\times B$, $r$ is a birational map and $f\colon \pp^1\to B/G$ is an
isomorphism. Again by proposition \ref{isostruttura}, the Albanese image of $V$ is
isomorphic to  $D/G$. So we have either $G=\{1\}$, corresponding to case i), or 
$2q(V)-1\le \gamma$,
with equality if and only if $G=\Z_2$ acts freely on $D$. The latter case corresponds to
case ii). \hfill\qed 

\section{General properties of generating pairs}\label{generalita}

In this section we give some useful information on the degree, genus and Kodaira dimension
of a generating pair. \par

\begin{nota} If $(h\colon V\to W, L)$ is a generating pair of degree $d$ and genus $g$, we
write $C$ for a general curve of $|L|$ and $C'=h^*C$, so that $C$ and
$C'$ are smooth curves of genera $g$ and $d(g-1)+1$ respectively, and
$h$ restricts to an unramified cover $\pi\colon C'\to C$ of degree $d$.\end{nota}

\begin{lem}\label{isolemma} Let  $(h\colon V\to W, L)$ be a generating pair of degree $d$
and genus $g$. If the Albanese image of $V$ is a curve, then
$d(g-1)+1>2q(V)-1$.
\end{lem} \proof According to proposition \ref{xiaose}, we distinguish three cases. 
Setting $D=C'$,
$\gamma=d(g-1)+1$ and keeping the rest of notation of proposition
\ref{xiaose}, we only need to exclude the occurrence of the first two
cases:\par

i) $V$ is ruled and $C'$ is a section: in this case the adjoint system $|K_V+L'|$ is
empty, contradicting assumption (GP3) of definition \ref{bpair};\par

ii) there are two subcases:
\begin{list} {ii-\alph{alpha})}{\usecounter{alpha}}
\item Assume $B=\pp^1$. Then $V$ is ruled over $C'/\Z_2$ and $C'$ is a
bisection of $V$ meeting each fibre of the map $p\colon V\to C'/\Z_2$ in two
distinct points interchanged by the free $\Z_2$ action. By repeatedly blowing
down $-1$ curves $E$ such that $EL'\le 1$, one obtains a map $f\colon  V\to
V'$ such that $V'$ is minimal and the map $p$ factors as $p'\circ f$, where
$p'\colon  V'\to C'/\Z_2$.  The curve $C''=f(C')$ is smooth and the induced map
$f\colon  C'\to C''$ is an isomorphism.
Moreover, the map $p'\colon  V'\to C'/\Z_2$ is a projective bundle, i.e.  there
exists a rank $2$  vector bundle $M$ on $C'/\Z_2$ such that
$V'=\Proj_{C'/\Z_2}(M)$, and $C''$ meets each fibre of $p'$ in two distinct
points interchanged by the free $\Z_2$ action.\par

If we denote by $H$ the tautological  section of $V'$ and by $L''$ the line bundle
determined by $C''$ on $V'$, then the condition that the projection map $C'\to
C'/\Z_2$ is unramified of degree $2$ is equivalent to  $L''$ being
numerically equivalent to $2H-\deg(M)F$, and thus we have  $L''^2=0$. This would imply
$L'^2\leq 0$, contradicting the fact that $L'$ is big.

\item Assume that $B$ is not rational. Notice that $(C'\times B)/\Z_2$ is the quotient of 
$C'\times B$ by a free $\Z_2$ action. Hence it is smooth. In addition it is
minimal, since it is a free quotient of the minimal surface 
$C'\times B$. This implies that the birational map $r\colon V\to (C'\times
B)/\Z_2$ is a
morphism. Let
$C''$ be a fibre of the morphism $(C'\times B)/\Z_2\to B/\Z_2=\pp^1$. Since, by
proposition \ref{xiaose}, $C'$ is the strict transform of $C''$ via $r$ and since
$C''^2=0$, we have again that $L'^2\leq 0$, which is impossible since $L'$ is
big.  \hfill\qed
\end{list}

\begin{lem}\label{qlemma} Let  $(h\colon V\to W, L)$ be  a generating  pair of degree $d$
and genus $g$. Then $d(g-1)+1\geq 2q(V)-1$, and if equality holds then the
Albanese image of $V$ is a surface. \end{lem} \proof Consider a pencil ${\cal
P}\subset |L|$ such that the general curve is  smooth and irreducible. Up to
blowing up, we may assume that the pull-back of ${\cal P}$ on $V$ via $h$ is a 
base point free pencil. If the  corresponding fibration is not isotrivial, then
the claim holds by proposition \ref{xiaofib}.  If the  fibration is
isotrivial,  then the Albanese image of $V$ is a
curve according to proposition \ref{isostruttura}, and by lemma
\ref{isolemma} we have
$d(g-1)+1> 2q(V)-1$. \hfill \qed

\begin{prop}\label{regular} If  $(h\colon V\to W, L)$ is  a generating pair of degree $d$,
then $q(W)=0$,  and the list of possibilities is as follows:
\begin{list} {\roman{rom})}{\usecounter{rom}} 
\item $d=2$, $q(V)=g-1$,
\item $d=3$, $g\le 3$, $q(V)=2(g-1)$
\item $d=4$, $g=2$, $q(V)=3$. 
\end{list}

\noindent If the pair is good, then case i) holds;  in case ii), $g=3$, and in case iii)
 the Albanese image of $V$ is a surface.
\end{prop} \proof By Kawamata-Viehweg's vanishing theorem one has  $h^0(W,
K_W+L)=\chi(W)+g-1$. Analogously, one has $h^0(V, K_V+L')=\chi(V)+d(g-1)$, and
thus $q(V)-q(W)=(d-1)(g-1)>0$ by condition (GP3) of definition \ref{bpair}.
Assume that $q(W)>0$. By proposition \ref{albcurva}, the Albanese image of $V$ is a curve
and lemma \ref{isolemma} implies that $d(g-1)+1> 2q(V)-1=2(d-1)(g-1)+2q(W)-1$, but this is
impossible, since $d, g\ge 2$. So, $q(W)=0$ and, according to lemma \ref{qlemma}, one has
$d(g-1)+1\geq 2q(V)-1=2(d-1)(g-1)-1$. The statement includes all
possible solutions. In cases ii) with $g=3$  and iii) we also apply lemma
\ref{qlemma}.\par 
Assume now that the pair is good.  By the above discussion, we
have $d\leq 3$. By (\cite{Beau}, Prop. 4.1 and Rem. 4.2), if $\psi\colon 
X\to S$ is a good canonical cover of degree $3$, then $q(X)\leq 3$. On the
other hand, by proposition \ref{explain}, canonical covers arising from a 
good generating pair of degree $3$ and
genus $g$ satisfy $q(X)=2(g-1)\ge 4$.\hfill\qed

\begin{prop}\label{kappa} Let  $(h\colon V\to W,L)$ be  a generating pair:
then $\kappa(V)=\kappa(W)$.\end{prop} \proof Remark first of all that
$\kappa(V)\ge \kappa(W)$. Hence we may assume $\kappa(W)\leq 1$.
Consider the following commutative diagram:
$$\begin{array}{rcccl}
\ &V'&\stackrel{b}{\rightarrow} &V&
\ \\
\scriptstyle{h'}\!\!\!\!\!\! & \downarrow &  &
\downarrow &
\!\!\scriptstyle{h}\\
\ &W' & \stackrel{f}{\rightarrow} & W &
\end{array}$$
\noindent where $f\colon  W'\to W$ is a minimal desingularization and $h'\colon  V'\to
W'$ is obtained by taking base change,  normalizing and finally solving the
singularities of the surface thus obtained.  We notice the following facts:\par

\noindent (i)  since $V$ and $V'$ are smooth surfaces, $b$ is a
sequence of blow-ups and thus $K_{V'}=b^*K_V+E$, where
$E$ is an effective divisor  supported on the $b$-exceptional locus. In addition,  for
every
$m\ge 1$  we have
$|mK_{V'}|=b^*|mK_V|+m E$. 
\noindent (ii) Since $W$ has only canonical singularities, one has 
$K_{W'}=f^*K_W$. Therefore we have $b^*K_V=b^*(h^*K_W)=h'^*K_{W'}$.  \par

Suppose that $\kappa(W)=-\infty$, i.e. $W$ is rational by proposition \ref{regular}.
Hence also $W'$ is rational, and therefore there is an effective irreducible big
divisor $D$ on $W'$ such that $D K_{W'}<0$. By remark (ii) above, there is an effective
big divisor $D'$ on $V'$ such that $D' (b^*K_V)<0$. This, together with  remark (i),
shows that $\kappa(V')=\kappa(V)=-\infty$. \par

Assume now that $\kappa(W)=0$; then there exists a nef and big line bundle $H$ on $W'$
such that $H K_{W'}=0$. Thus $(h'^*H) (b^*K_{V})=(h'^*H)(h'^*K_{W'})=0$, 
and thus $h'^*H$ is a nef and big divisor that has zero intersection with the moving part
of any pluricanonical system.  Thus it follows  that $\kappa(V')\le 0$. If $\kappa(W')=1$,
then there exists a fibration $f\colon W'\to D$, where $D$ is a smooth curve, such that
the general fibre $E$  of $f$ is an elliptic curve. So $(h'^*E) (b^*K_V)=0$, and thus the
maps given by the pluricanonical systems are all composed
with the fibration $f'=f\circ h$. This shows that
$\kappa(V)\le 1$. \hfill\qed

\medskip
According to the previous proposition, we may, and will, speak of the {\it
Kodaira dimension} of a generating pair. Generating pairs of degree $2$ and
Kodaira dimension $0$ are completely described in proposition
\ref{kappa0}.

\section{Pairs of degree  $2$ and Prym varieties}\label{grado2}

We consider here the case of generating pairs of degree $2$. The relevance of this case is
underlined in proposition \ref{regular}, where it is shown that all good
generating pairs have degree $2$ and that all
generating pairs of degree $>2$ have genus $\le 3$.\par

If $C\in |L|$ is a general curve and $C'=h^*C$, then the map $h$ induces an
\'etale double cover $\pi\colon C'\to C$. If one denotes by $J$  (resp.
$J'$) the Jacobian of $C$ (resp. $C'$), then the connected  component  of the kernel of the norm map $\pi_*\colon J'\to J$  is
a $(g-1)$--dimensional containing the
origin is an abelian variety,  on which  the principal
polarization of $J'$ induces  the double of a principal polarization. This
principally polarized abelian variety is called the {\em Prym variety} of
$C'\to C$ and it is denoted by $Prym (C',C)$. The connection between generating
pairs and Prym varieties  is explained in the following theorem.

\begin{thm}\label{prym} Let $(h\colon V\to W,L)$ be a generating pair of degree 2. Let
$C\in |L|$ be a general curve. Then there is a natural isomorphism $\varphi\colon 
Prym(C',C)\to A$, where $A=Alb(V)$ is the Albanese variety of $V$. In particular
$Prym (C',C)$ does not depend on $C\in |L|$. 
\end{thm} \proof  Under the present assumption, the singular points of $W$
form a set of $t$ ordinary double points, where $t$ satisfies the relation
$\chi(V,\OO_V)=2 \chi(W,\OO_W)-t/4$.  Evaluating the Euler characteristic of
$V$ and $W$ as in proposition \ref{regular}, one deduces that
$t=4(g+p_g(W))>0$. So, one can choose a ramification point $x_0$ for $h$ in
$V$. Since $W$ is regular by proposition \ref{regular}, the Albanese map of
$W$ with base point $h(x_0)$ is the zero map. The Albanese map of $V$ with
base point $x_0$, denoted by $\alpha\colon V\to A$, is equivariant with
respect to the involution induced by $h$ and the multiplication by $-1$ on
$A$. In particular, the restriction $\alpha|_{C'}\colon C'\to A$ is also
equivariant.\par

Now we use the universal property of Prym varieties (cf. \cite{LB}, page
382). Let $\beta\colon C'\to Prym(C',C)$ be the Abel--Prym map with respect to
a point $c'\in C'$  and let $\tau\colon A\to A$ be the translation by
$\alpha(c')$. Then there is a unique homomorphism $\varphi\colon Prym(C',C)\to
A$, independent of $c'\in C'$, such that $\alpha|_{C'}=\tau\circ \varphi\circ
\beta$.\par

Denote by $J'$ the Jacobian of $C'$. Let $j\colon  C'\to J'$ be the Abel map
with base point $c'$ and $\gamma\colon  J'\to Prym(C',C)$ the map such that
$\beta=\gamma\circ j$. Let $i_*\colon J'\to A$ be the homomorphism induced by
the inclusion $i\colon C'\to V$ and the choice of $c'\in C'$. Notice that, up
to a translation, we have $\alpha_{|C'}=i_*\circ j$. Then it is clear that
$i_*$ factors, up to a translation, as $\varphi\circ \gamma$. The
differential of $i_*$ at the origin of $J'$ is dual to the map $H^1(V,
\OO_V)\to H^1(C', \OO_{C'})$, which is injective since $H^1(V, \OO_V(-L'))=0$
because $L'$ is big and nef. So $i_*$ is surjective and $\varphi$ is an
isogeny since $A$ and $Prym(C',C)$ both have dimension $g-1$ by proposition
\ref{regular}. To show that $\varphi$ is an isomorphism, it is enough to prove
that $i_*$ has connected fibres. In turn, this follows if we show that the map
$H_1(C', \Z)\to H_1(V,\Z)$ induced by the inclusion $i\colon C'\to V$ is
surjective. The system  $|L'|$ has no fixed part by assumption, so by theorem
$6.2$ of \cite{zariski} there exists an integer $k$ such that $|kL|$ gives a
morphism $g:V\to \pp^N$; the image of $g$ is a surface, since $L'$ is big. So
there exists an hyperplane $H$ in $\pp^N$ such that $g\inv H=C'$ {\em as
sets}. By Theorem $1.1$, page $150$, of \cite{GM}, the map
$\pi_1(C)\to\pi_1(V)$ is surjective, and thus  $H_1(C',\Z)\to H_1(V,\Z)$ is
surjective too. \hfill\qed 

\begin{cor}\label{albsurface} Let $(h\colon V\to W,L)$ be a  generating
pair of degree $2$ and  genus $g$; then the Albanese image of $V$ is a
surface. In particular, the Kodaira dimension of the pair is non--negative. 
\end{cor} \proof Assume that the Albanese image of $V$ is a curve $\Gamma$.
Then $\Gamma$ has genus $g-1$. On the other hand, by theorem \ref{prym}, the
Albanese image of $V$ contains the Abel--Prym image  of $C'$, which is
isomorphic to $C'$ (cf. \cite{LB}, prop. 12.5.2), since $C'$ is not hyperelliptic. This
is a contradiction and thus the claim is proven. 
\hfill\qed

\begin{cor}\label{mintame} Let $(h\colon V\to W, L)$ be a    generating
pair of degree $2$; then $(h\colon V\to W, L)$ is obtained from a minimal pair
by a sequence of simple blow-ups of weight $0$ or $1$.
\end{cor}\proof  Denote by $i\colon V \to V$ the involution induced by $h$ and  let  $E$ be 
a $-1$ curve of $V$. We claim that either $L'E=0$ or $L'E=1$. Let $\epsi\colon
V\to V_0$ be the blow--down of $E$, let $C'\in |L'|$ be smooth and let
$C_0=\epsi(C')$;  notice that  $C_0$ is singular if and only if $L'E>1$. Let
$\alpha\colon V\to A$ be the Albanese map of $V$; $A$ is also the Albanese
variety of $V_0$ and, if we denote by $\alpha_0\colon V_0\to A$ the Albanese
map of $V_0$, one has $\alpha=\alpha_0\circ\epsi$. Thus 
$\alpha(C')=\alpha_0(C_0)$; by theorem \ref{prym}, $\alpha(C')$ is isomorphic
to $C'$,  since $C'$ is not hyperelliptic, and thus $C_0$ is smooth and
$L'E\le 1$. Let $E'$ be the image of $E$ via $i$; $E'$ is also a $-1$--curve
and thus, since $\kappa(V)\ge 0$ by corollary \ref{albsurface}, either $E=E'$
or $E$ and $E'$ are disjoint. If $E=E'$, then $E$ contains precisely $2$ fixed
points of $i$, but this contradicts the fact that $E^2$ is odd. So $E\ne E'$
and $F=h(E)=h(E')$ is a $-1$ curve contained in the smooth part of $W$. Let
$V'$ be the surface obtained by blowing down  $E$ and $E'$, let $W'$ be the
surface obtained by blowing down $F$ and let $h'\colon V'\to W'$ be the double
cover induced by $h$; if one denotes by $M$ the direct image of $L$, then it
is easy to check that $(h'\colon V'\to W', M)$ is also a   generating pair. By
iterating this process finitely many times, one eventually obtains a 
generating pair with $V$ minimal. Thus $K_V=h^*K_W$ is nef, and it follows
that $K_W$ is also nef and $W$ is minimal, too. \hfill\qed

\begin{cor}\label{pg>0} Let $(h\colon V\to W,L)$ be a  generating pair of genus
$g$ and degree $2$. Then:\par

(i) $p_g(V)=p_g(W)\ge g-2>0$;\par 
(ii) if the Kodaira dimension of the pair is $2$, then $p_g(V)=p_g(W)\ge max
\{g-1,2g-6\}$; if $p_g(V) =2g-6$ then $V$ is birational to the product of a curve of
genus $2$ and a curve of genus $g-3$.
\end{cor} \proof By corollary \ref{albsurface}, $g-1=q(V)>1$. Thus we have $\chi(V)\ge 0$,
$p_g(W)=p_g(V)\ge q(V)-1=g-2>0$. The case of Kodaira dimension $2$ follows from the
theorem at pg. 345 of \cite {Beau4}. \hfill\qed 

\begin{cor}\label{moduli}  Let $(h\colon V\to W,L)$ be a  generating pair of genus
$g$ and degree $2$. If $|C'|\subset
|L'|$ is  a pencil containing a smooth curve, then $|C'|$ is not isotrivial.
\end{cor} \proof Follows from corollary \ref{albsurface} and proposition
\ref{isostruttura}. \hfill\qed

\medskip
If we denote by ${\cal R}_g$ the moduli space of \'etale double covers of
curves of genus $g$  and by ${\cal A}_{g-1}$ the moduli space of principally
polarized  abelian varieties of dimension $g-1$, then the {\em Prym map}
${\cal  P}_g\colon {\cal R}_g\to {\cal A}_{g-1}$ associates to every
isomorphism class of \'etale double covers the corresponding Prym variety.
The  geometry of Prym varieties has been extensively studied  by many authors.
We are going to use some of these results in order to give a bound on the
genus of good generating pairs. 

\begin{prop} \label{notfinite} Let $(h\colon V\to W,L)$ be a  generating
pair of genus $g$ and degree $2$. Let $C\in|L|$ be general and let $C'=h^*C$.
Then the fibre of the Prym map ${\cal P}_g\colon {\cal R}_g\to {\cal A}_{g-1}$
at the point of ${\cal R}_g$  corresponding to the double cover $C'\to C$ has
positive dimension. \end{prop} \proof Follows from theorem \ref{prym} and
corollary \ref{moduli}. \hfill\qed

\medskip
It is known that the Prym map is generically finite for
$g\ge 6$ (cf. the survey \cite{Beausurvey} and the references quoted
therein). However there exist positive dimensional fibres of ${\cal P}_g$ for
any value of $g$. In order to state Naranjo's theorem \ref{fibreprym} that
characterizes the positive dimensional fibres of ${\cal P}_g$ for high values
of $g$, we recall that a curve $C$ is called bi-elliptic if and only if it
admits a double cover $C\to E$ onto an elliptic curve $E$.

\begin{thm}\label{fibreprym} (Naranjo, see \cite{Na2}, page 224 and \cite
{Na1}, theorem (10.10)) Let
$C'\to C$ be an unramified double cover of a genus $g$ curve $C$.\par

(i) If $g\geq 13$, then the fibre of ${\cal P}_g$
at the point of ${\cal R}_g$ corresponding to   $C'\to C$ is positive
dimensional if and only if  $C$ is either hyperelliptic or  is bi-elliptic.
In addition, in the latter case, if $C\to E$ is a double cover of an elliptic
curve, then the Galois group of the composition $C'\to C\to E$ is 
$G=\Z_2\times \Z_2$ and each quotient of $C'$ under an element of $G$ has
genus strictly greater than $1$.\par

(ii) If $g\geq 10$,  the fibre of ${\cal P}_g$ at the point of ${\cal
R}_g$ corresponding to $C'\to C$ is positive dimensional and $C$ is
bi-elliptic, then the Galois group of the composition $C'\to C\to E$ is 
$G=\Z_2\times \Z_2$, and each quotient of $C'$ under an element of $G$ has
genus strictly greater than $1$. 
\end{thm} 

From the point of view of generating pairs, the hyperelliptic case in Theorem
\ref{fibreprym} corresponds to the case of generating pairs of degree $2$ which are not
good, and example \ref{product} shows that these exist for arbitrary values of $g$. On
the other hand, the bielliptic case can be excluded for good generating pairs with $g$
large, as theorem \ref{g<12} below shows. 

We recall some general and elementary properties of bi-elliptic curves and bi-double
covers, i.e. finite covers with Galois group $\Z_2\times \Z_2$ (cf. \cite{Na1}, page 50
and ff.; \cite {ritaabel}). If $C$ is bi-elliptic, then the double cover $C\to E$ with
$E$ elliptic is unique up to automorphisms of $E$ if $g\ge 6$. Analogously, a bi-elliptic
curve $C$ is not hyperelliptic if $g\ge4$ and it is not trigonal if $g\ge 6$.\par

If $C'\to C$ is an \'etale double cover of a bi-elliptic curve $C\to E$, then  the
composition $C'\to C\to E$ is a degree $4$ cover of $E$ whose Galois group $G$ contains
$\Z_2$. Assume that $G=\Z_2\times\Z_2$,  denote by $\sigma$ the element of  $G$ such that
$C'/<\sigma>=C$ and by $\sigma_i$ ($i=1,2$) the remaining non trivial elements. For
$i=1,2$, set $p_i\colon C'\to C_i=C'/<\sigma_i>$ the corresponding projection and notice
that $C_i$ is a smooth curve of genus $g_i$, where $g_1+g_2=g+1$. Then there exists a
cartesian diagram:  
\begin{equation}\label{biell}
\begin{array}{rcccl}
\phantom{1} &C'&\stackrel{\pi_2}{\rightarrow} &C_2&
\phantom{1} \\
\scriptstyle{\pi_1}\!\!\!\!\!\! & \downarrow & \phantom{1} &
\downarrow &
\!\!\!\!\!\!\scriptstyle{\phi_2}
\\
\phantom{1} & C_1 & \stackrel{\phi_1}{\rightarrow} &E &
\phantom{1}
\end{array}
\end{equation} where, for $i=1,2$, $\phi_i\colon C_i\to E$ is a double cover,
and the branch loci $\Delta_i$ of $\phi_i$, $i=1,2$, are disjoint. Moreover,
$\sigma=\sigma_1\circ \sigma_2$. The group $G$ also acts on $Prym(C',C)$. We
denote  by $P_i$ the connected component containing the origin of the fixed
locus of the action of $\sigma_i$ on $Prym(C',C)$ ($i=1,2$), and we observe
that $(P_1, P_2)$ is a pair of complementary abelian subvarieties  of
$Prym(C',C)$ of dimensions $g_1-1$ and $g_2-1$, respectively.

\begin{lem}\label{bielliptic} Let $(h\colon V\to W,L)$ be a good generating
pair of genus $g\geq 10$. Then the general curve $C\in |L|$ is not bielliptic.
\end {lem} \proof By corollary \ref{mintame} we may assume that the pair is
minimal. Suppose, by contradiction, that the general curve $C\in |L|$ admits
an elliptic involution $C\to E$, which, as we saw, is unique up to
automorphisms of $E$. Moreover, by part (ii) of theorem \ref {fibreprym} and
by proposition \ref{notfinite}, the Galois group $G$ of the composition $C'\to
C\to E$ can be identified with $\Z_2\times \Z_2$. Theorem \ref {fibreprym}
also ensures that there exists a cartesian diagram as in (\ref{biell}), with
$C_i$\, of genus $g_i>1$. We wish to extend this construction to $V$. \par

In order to do this, we  prove first that we may choose the involutions
$\{\sigma_1, \sigma_2\}$ consistently on the curves $C'=h^*C$ as $C$ varies in
$|C|$. In other words, there is a double cover $\Psi\to \Phi$ of the open
subset $\Phi$ of $|C|$ parametrizing smooth curves, such that its fibre at a
general point $C\in\Phi$ is the pair of involutions $\{\sigma_1, \sigma_2\}$
acting on $C'=h^*C$. We want to prove that $\Psi$ is the union  of two
irreducible components both mapping birationally to $\Phi$. In order to do
this, we have to prove that there are two sections of $\Phi\to \Psi$ mapping
the general point $C\in \Phi$ to $\sigma_1$, resp. $\sigma_2$, namely that we
can rationally distinguish $\sigma_1$ from $\sigma_2$.\par

Recall that, by theorem \ref{prym}, $Prym(C',C)$ is isomorphic, in a
canonical way, to the Albanese variety $A$ of $V$. In this isomorphism, the
connected component $P_i$ of the origin of the fixed locus of the action of
$\sigma_i$ on $Prym(C',C)$  maps to an abelian subvariety $B_i$ of $A$
($i=1,2$).  The pair $(B_1,B_2)$ of complementary subvarieties can vary only
in a discrete set, and therefore it is constant, independent of $C$. This
proves our claim about the reducibility of $\Psi$. \par

Next we claim that there are  involutions $\tau_i$ on $V$ inducing $\sigma_i$
on the general $C'$, for $i=1,2$. Indeed, let ${\cal F}$ be a general pencil
inside $|C|$. If $x\in V$ is a general point, define $\tau_i(x)$ as
$\sigma_i(x)$, where $\sigma_i$ is the involution defined on the unique curve
$C'$ in $h^*({\cal F})$ passing through $x$. Since $V$ is minimal, $\tau_i$
extends to an automorphism of $V$. Notice that $\tau_i$ is independent of
${\cal F}$, otherwise,  as ${\cal F}$ varies in a general rational
$1$--parameter family of pencils, the point  $\tau_i(x)$, $x\in V$ general, 
would describe a rational curve, hence $\kappa(V)$ would be negative, against
proposition \ref{pg>0}.\par We denote by $S_i$ the quotient surface
$V/<\tau_i>$, by $h_i\colon V\to S_i$ the projection onto the quotient and by
$C_i$ the image in $S_i$ of a general $C'$. The singularities of $S_1$ and
$S_2$, if any, are $A_1$ points and $q(S_i)=g_i-1$. By proposition
\ref{xiaose}, if the curves $C_i$ vary in moduli, then  $g_i\le 3$, thus
$g=g_1+g_2+1\le 7$, a contradiction. If the curves $C_i$ do not vary in
moduli, then the Albanese image of $S_i$ is a curve by proposition
\ref{isostruttura} and  the  inequality $g_i\le 3$ ($i=1,2$) holds by
proposition \ref{xiaofib}, since $q(S_i)\neq g_i$. \hfill\qed

\medskip

Now we are ready to prove the following basic result:

\begin{thm}\label{g<12} Let $(h\colon V\to W,L)$ be a good generating pair of genus $g$.
Then $g\le 12$, $q(V)\le 11$. 
\end{thm} \proof Suppose, by contradiction, that $g\geq 13$. According to proposition
\ref{notfinite} and to part (i) of theorem
\ref{fibreprym}, we can assume that the general $C\in
|L|$ is bi-elliptic. This, on the other hand, contradicts lemma
\ref{bielliptic}. \hfill\qed\medskip

\begin{cor}\label{K2}  Let $(h\colon V\to W,L)$ be a good generating pair of genus $g$.
Assume that, in addition, $V$ is of general type. Then $K_W^2\le 529$. 
\end{cor} \proof Follows by applying the index theorem to $L$ and $K_W$ on $W$. \hfill\qed

\medskip

A more precise statement is the following:

\begin{prop}\label{genere} Let $(h\colon V\to W,L)$ be a good generating pair of genus
$g$. Then $g\le 12$ and:
\begin{list}{(\roman{arab})}{\usecounter{arab}}
\item if the general curve $C$ in $|L|$ is bi-elliptic or trigonal, then $g\le
9$ and $q(V)\le 8$;
\item if $10\le g\le 12$ then either the general curve $C$ in $|L|$ is a
smooth plane sextic (and $g=10$) or it is not bi-elliptic and has a
base point free $g^1_4$.
\end{list}
\end{prop} \proof The proof follows from theorem \ref{g<12}, theorem \ref{fibreprym}, and
from the following results: 

\begin{list}  {(\roman{arab})}{\usecounter{arab}}
\item (Green-Lazarsfeld \cite{GL}) Assume $g\geq 10$. If the fibre of the Prym
map ${\cal P}_g$ is positive dimensional at the point of ${\cal R}_g$ 
corresponding to a double cover $C'\to C$, then either $C$ has
a $g^1_4$ or it is a smooth plane curve of degree six (and genus 10).
\item (Naranjo \cite{Na2}) Assume $g\ge 10$. Then the fibre of ${\cal P}_g$
over the point of ${\cal R}_g$ corresponding to a double cover of a trigonal
curve $C$ is finite. \hfill\qed
\end{list}

\section{Good generating pairs with $h^0(W,L)\ge 3$}\label{sezgood} 

This section is devoted to the proof of the following:
\begin{thm}\label{maingood} Let $(h\colon V\to W, L)$ be a good generating
pair such that $L$ is ample and $h^0(W,L)\ge 3$; then $(h\colon V\to W, L)$ is
one of the following:\par

\noindent
(i) example \ref{Beauville}, and in this case 
$h^0(W,L)=4$, $L^2=4$, $g=3$;\par
\noindent
(ii) a blow--up of weight $1$ of case (i) above, and in this case
$h^0(W,L)=3$, $L^2=3$, $g=3$;\par

\noindent
(iii) example \ref{Rita}, and in this case $h^0(W,L)=3$, $L^2=3$, $g=4$.
\end{thm} 
Theorem \ref{maingood} will follow from a series of auxiliary results
(proposition \ref{classificagood} and theorem \ref{trigonalcase}), containing
also some additional information on generating pairs. We also make use of the
following result, which is proven in section \ref{altricasi} (it
follows from proposition \ref{kappa0} and \ref{kappa1}).

\begin{prop}\label{kodairagood} If $(h\colon V\to W)$ is a good generating
pair with $\kappa(V)\le 1$ and $h^0(W,L)\ge 3$, then it is obtained from

example \ref{Beauville} by a sequence of blow-ups, at most one of which is
essential, of weight $1$. \end{prop}

We start by using  Reider's method to give an upper bound for $L^2$ for most
generating pairs.
\begin{prop}\label{Reider} If $(h\colon V\to W, L)$ is a generating pair  of
degree $2$ and non--negative  Kodaira dimension, then $L^2\leq 4$. \end{prop}
\proof Since $L'^2=2 L^2$, it suffices to show that $L'^2\leq 9$. We assume
that $L'^2\geq 10$ and we observe that, by the  hypothesis, the linear system
$|K_V+L'|$ is not birational on $V$. Indeed, if $x\in W$ is a general point
and $h^{-1}(x)=\{x_1, x_2\}$, then $x_1$, $x_2$ are identified by $|K_V+L'|$.
According to Reider's Theorem (\cite{Reider}, Thm. 1 and Cor. 2), there exists
an effective divisor $B=B_x$ passing through $x_1$ and $x_2$, such that $L'
B=1$ or $2$ and $B^2\le 0$. Since $x$ is general we must have $B^2=0$ and, by
standard arguments, we may assume that $B$ moves in a base point free pencil
and $L' D>0$ for each component of a general $B$ . Since the general curve
$C'$ of $|L'|$ is irreducible and meets the general curve $B_x$ at the  points
$x_1$ and $x_2$, it follows that $L' B=2$. If the general $B$ is reducible,
then $B=B_1+B_2$, where $B_1$, $B_2$ are numerically equivalent irreducible
curves. Then $L' B_i=1$ and $V$ is covered by rational curves, contradicting
the assumption $\kappa(V)\geq 0$. So, $L' B=2$ and  $B$ is irreducible.
Furthermore, the general fibre of $h$ is contained in some curve of the pencil
described by $B_x$, as $x$ varies in $W$. This immediately implies that each
curve of this pencil is invariant under the involution $\iota$ determined by
$h$. On the other hand $|L'|$ cuts  on a general curve $B$ a $g_2^1$, which of
course induces on $B$  the restriction of $\iota$. This means that the image
of $B$ via $h$ is a rational curve on $W$, which therefore has a pencil of
rational curves. But this is a contradiction to $\kappa(W)\geq 0$. \hfill\qed

\begin{lem}\label{classifica} Let $(h\colon  V\to W,L)$ be a  generating pair
such that $h^0(W,L)\geq 3$. Then there are the following possibilities:\par
(i) $h^0(W, L)=3$ and $2\leq L^2\leq 4$;\par

(ii) $h^0(W,L)=4$, $L^2=4$ and $|L|$ is base point free.
\end{lem}  \proof If $h^0(W, L)=r$, then the restriction of the system $|L|$
to a general $C$ is a linear system $|D|$ of dimension $r-2>0$ and  degree
$L^2\leq 4$, according to proposition \ref{Reider}. We denote by $|M|$ the
moving part of $D$.  If $r=3$, then $L^2\ge \deg M\ge 2$, since $C$ is not
rational, and (i) is proven. If $r>3$, then $4\ge  L^2\ge \deg M\geq
2\dim|M|=2(r- 2)\ge 4$. Thus $L^2=\deg M=r=4$ and $|L|$ is base point free.
\hfill\qed

\begin{prop}\label{classificagood} Let $(h\colon  V\to W,L)$ be a good
generating pair such that $h^0(W,L)\geq 3$. The possible cases are:  
\begin{list} {(\roman{arab})}{\usecounter{arab}} 
\item $L^2=4$, $h^0(W,L) = 3$ and $|L|$ is base point free. In particular,
the general $C\in |L|$ is tetragonal.
\item  $h^0(W,L) = 3$ and either $L^2=3$ or $L^2=4$ and $|L|$ has a simple
base point. In particular, the general $C\in |L|$ is trigonal;
\item  $L^2=4$, $h^0(W,L) = 4$, and the pair is obtained from Beauville's
example \ref{Beauville} via unessential blow-ups. \end{list}
\end{prop} \proof We denote by $|D|$ the restriction of $|L|$ to a general
$C$ of $|L|$. Assume that we are in case (i) of lemma \ref{classifica}: then
(i) and (ii) follow  by remarking that the moving part of $|D|$ has degree
$>2$, since $C$ is not hyperelliptic.\par

Assume that we are in case (ii) of lemma \ref{classifica}. Then Clifford's
theorem implies that $g=3$, $|D|$ is the canonical system and $C$ is embedded
by $|D|$ as a smooth plane quartic. So the linear system $|L|$ maps the
surface $W$ birationally onto a quartic $Q\subset \pp^3$. The Kodaira
dimension of $W$ is non--negative by corollary \ref{pg>0} and thus it is zero.
Claim (iii) now follows by proposition \ref{kodairagood}. \hfill\qed 

\begin{prop}\label{quattro} 
If $(h\colon V\to W, L)$ is a good generating pair of Kodaira dimension 2 such
that $L^2=4$ and $h^0(W,L)\le 3$, then  $10\leq g\leq 12$.
\end{prop} \proof Notice first of all that the inequality $g\le 12$ follows
from theorem \ref{g<12}.\par

By corollary \ref{mintame} it follows
that the pair is obtained from a minimal pair by unessential blow--ups. Thus
we may assume that the pair is minimal. Write $h^0(W,L)=2+l$, so that either
$l=0$ or $l=1$.\par

Since $W$ is of general type, one  has $0<K_W L=2g-2-L^2=2g-6$, hence $g\ge 4$. For a
general $C\in |L|$, consider the exact sequence:\par

$$0\to H^0(W,K_W-L)\to H^0(W,K_W)\to H^0(C,K_C-L_{|C})$$

\noindent  and notice that $h^0(C,K_C-L_{|C})=g-4+l$, by the regularity of
$W$ and by Riemann--Roch applied to $C$.\par

Assume first $g\le 6$. Then by corollary \ref {pg>0}, we have
$h^0(W,K_W-L)\geq g-1-(g-4+l)= 3-l\ge 2$. Notice that $h^0(W,K_W-2L)=0$, since
$(K_W-2L)L=2g-14< 0$ and $L$ is nef. Thus $|K_W-L|$ cuts out  on $C$ a  linear
series of  dimension $\ge 1$ and of degree $2g-10\le 2$, contradicting the
assumption  that $C$ is not hyperelliptic.\par

Therefore we have $g\ge 7$. Corollary \ref{pg>0} and the above exact sequence yield: 
$h^0(W,K_W-L)\ge 2g-6-(g-4+l)= g-2-l$. Let $h^0(W,K_W-2L)=r$. Then $|K_W-L|$ cuts out on
$C$ a special linear series of dimension $g-3-r-l$ and degree $2g-10$. By Clifford's
theorem and by the fact that $C$ is not hyperelliptic, we  have
$g-3-r-l<g-5$, namely  $r+l\ge 3$. Thus we either have $l=0$, $r\ge 3$ or $l=1$, $r\ge 2$.
If $g\le 8$, then $L(K_W-3L)=2g-18< 0$, and thus $h^0(W, K_W-3L)=0$ and
$|K_W-2L|$ cuts out on $C$ a special linear series of dimension at least $r-1\geq 1$ and
of degree $2g-14\le 2$, contradicting again the assumption that $C$ is not hyperelliptic. 
\par

If $g=9$, then $L(K-3L)=0$  and thus $h^0(W, K_W-3L)\le 1$, since $L$ is nef and big.
Assume that $h^0(W,K_W-3L)=1$;  then we have $K^2_W-36=K_W(K_W-3L)\ge 0$, since $W$ is of
general type. On the other hand, the index theorem gives $K_W^2\le 36$. It follows that
$K^2_W=36$ and $K\sim_{num} 3L$. Therefore $K_W=3L$, since $K_W-3L$ is effective.  So
$r=h^0(W, K_W-2L)=h^0(W,L)=2+l$, and for $l=0$ this  contradicts the  above inequality
$r+l\ge 3$. If $l=1$, consider the exact sequence:
\begin{equation}\label{seq} 0\to (k-1)L\to kL\to kL_{|C}\to 0.\end{equation} 
By Clifford's theorem we have $h^0(C, 2L_{|C})\le 4$. So, for $k=2$, the  sequence
(\ref{seq}) implies $h^0(W,2L)\le 7$. Using this and sequence (\ref{seq}) for $k=3$, one
gets $p_g(W)=h^0(W,3L)\le 13$, and thus $\chi(V)=2+p_g(W)-g\le 6$. On the other hand,
Miyaoka--Yau's inequality would give $72=K^2_V\le 9\chi(V)=54$, a
contradiction.\par 

So we are left with the case $h^0(W, K_W-3L)=0$. If $l=0$, then $r\ge 3$ and  the
restriction of  $|K_W-2L|$ to $C$ is a $g^2_4$, contradicting again the fact that $C$ is
non-hyperelliptic of genus $9$. Thus the case $g=9$ and $l=0$ does not occur.\par

If $l=1$, then we have $r=2$, since for $r>2$ we can argue as above and show
that  $|K_W-2L|$ restricts to a $g^2_4$ on $C$. So we have $h^0(W,K_W-L)\le
2+h^0(C, (K_W-L)_{|C})\le 2+4=6$, where the last inequality follows again by
Clifford's theorem. In turn, $p_g(W)\le 6+h^0(C, {K_W}_{|C})=12$. On the other
hand, by proposition \ref{pg>0}, one has $p_g(W)\ge 12$, with equality holding
iff $V=C_1\times C_2$, with $C_1$ a curve of genus $2$ and $C_2$ a curve of
genus $6$. Thus $p_g(W)=12$ and the restriction map $H^0(W, K_W)\to
H^0(C,{K_W}_{|C})$ is surjective.   Since the canonical map of $V$ factors
through the map $h\colon V\to W$ and $q(W)=0$, the curve $C_2$ is
hyperelliptic,  and  the canonical map of $V$ has degree $4$.  The canonical
image $\Sigma$ of $V$ (and of $W$) is  $\pp^1\times \pp^1$ embedded via the
system $|\OO_{\pp^1\times\pp^1}(1,5)|$. The curves of $|L|$ are mapped
$2$--to--$1$ onto curves $D$  of $\Sigma$ with $D^2=2$. So the curves $D$ are
of type $(1,1)$, and thus rational. It follows that the general curve  of
$|L|$ is hyperelliptic, contradicting again the assumption that the pair be
good. \hfill\qed

\medskip
\begin{cor}\label{quattrocor}
If $(h:V\to W,L)$ is a good generating pair such that $L^2=4$ and
$h^0(W,L)=3$, then $10\le g\le 12$.
\end{cor}  \proof
Proposition \ref{kodairagood} implies that the Kodaira dimension of the pair
is $2$. The thesis follows from proposition \ref{quattro}.
\hfill\qed

\medskip

The next result shows that case (i) of
proposition \ref{classificagood} does not occur for $L$ ample.

\begin{prop}\label{basepoint}Let $(h\colon V\to W, L)$ be a good generating pair with
$h^0(W,L)=3$, $L^2=4$ and $L$ ample; then $|L|$ has one simple base point.
\end{prop} \proof  First of all we remark that the assumption that $L$ is
ample, corollary \ref{mintame} and proposition \ref{Reider} imply that the
pair is minimal. By lemma \ref{classificagood}  we only have to exclude that
$|L|$ is base point free.  So we  assume that $|L|$ has no base points and we
show that this leads to a contradiction. As usual we denote by $C$ a general
curve of $|L|$ and by $C'$ the inverse image of $C$ via $h$; we denote by
$\phi\colon W\to\pp^2$ the finite degree $4$ morphism given by $|L|$. Notice
that $\phi$ is flat, since it is a projective morphism with finite fibres from
a normal surface to a smooth one. Our proof requires various steps.\par

\noindent{\em Step 1: the polarized Abelian variety $Prym(C',C)$ is not a
Jacobian or  a product of Jacobians.} By theorem (4.10) of \cite{Beau0} and the
assumption that $C$ is not hyperelliptic, if $Prym(C',C)$ is a Jacobian then
one of the following holds: (i) $C$ is trigonal, (ii) $C$ is bielliptic, (iii)
$g\leq 6$. Case (iii) is impossible by corollary \ref{quattro}. Since $C$ has a free
$g^1_4$, case (i) implies $g\leq 6$, and therefore it is also excluded. Finally
case (ii) is excluded by lemma \ref{bielliptic}. \par

\noindent {\it Step 2: The curves of $|L'|$ are  $2$-connected.}
Notice first of all that the curves of $|L'|$ are $1$-connected since $|L'|$
is ample. Assume that $D\in |L'|$ is not $2$-connected, namely that $D=A+B$
with $A$, $B$ effective and $A B=1$. Then $1\leq L'A=A^2+AB=A^2+1$ and $1\leq
L'B=B^2+1$, hence $A^2\geq 0$, $B^2\geq 0$. If, say, $A^2=0$, then $L'A=1$ and
$A$ is irreducible, since $L'$ is ample, and rational, since $|L|$ is base
point free. This contradicts $\kappa(V)=2$, hence $A^2$ and $B^2$ are both
positive. Since $8=L'^2=A^2+B^2+2$, we have $A^2+B^2=6$ whereas the index
theorem gives $A^2B^2\leq (AB)^2=1$.\par

\noindent {\it Step 3: The branch divisor $Z$ of $\phi$ in ${\pp^2}$ is not a
union of lines.}
Here we  need to consider the intersection number of Weil divisors on $W$. We
recall that, since the singularities of $W$ are $A_1$ points, given  Weil
divisors $A$, $B$ on $W$, the intersection number $AB$ is an element of
$\frac{1}{2}\Z$, and it is an integer whenever $A$ or $B$ is Cartier.  Assume
that $Z$ is a union of lines and  let $R$ be a line contained in $Z$. Then
$C_0=\phi^*(R)\in |L|$ is of the form $C_0=mA+B$, with $2\leq m\leq 4$ and
with $A$, $B$ effective, non-zero Weil divisors such that $A$ is irreducible
and not contained in $B$. We set $C'_0=h^*C_0$, $A'=h^*A$, $B'=h^*B$, so that
$C'_0=mA'+B'$.\par
Notice that $4=L^2=mLA+LB\geq 2LA$ yields $1\leq LA\leq 2$. 
Assume first $LA=2$. Then one has  $B=0$, $m=2$, and thus $C'_0=2A'$. Recall
that by proposition \ref{prym} the abelian variety $Prym(C',C)$ is naturally
isomorphic to the Albanese variety of $V$ and denote by $\alpha\colon V\to
Prym(C',C)$ the Albanese map. If $\Xi$ is the principal polarization of
$Prym(C',C)$, then by Welters criterion $\alpha_*C'$ is homologically
equivalent to $\frac{2}{(g-2)!}\wedge^{g-2}\Xi$ for every $C'\in|L|$.
Thus $\alpha_*A'$ is homologically equivalent to
$\frac{1}{(g-2)!}\wedge^{g-2}\Xi$, and it follows that $A'$ is smooth and 
$Prym(C',C)$ is isomorphic to the Jacobian of $A'$. This is impossible by step
1, and therefore $LA=1$. The condition $LA=1$ implies that $m\le 3$, $B$ is
nonempty and $A$ is smooth and irreducible. Assume that
$m=3$, and let $R_1\subset Z$ be another line;  write $\phi^*R_1=m_1A_1+B_1$, with $2\le
m_1\le 3$, $A_1$ irreducible and not contained in $B_1$. The equality 
$1=A_1L=3A_1A+A_1B$ gives $AA_1=0$, $A_1B=1$, and thus $1=BL=mBA_1+BB_1\ge m\ge 2$, a
contradiction.  Thus, for every line $R\subset Z$, we have
$\phi^*R=2A+B$, with $A$ irreducible and not contained in $B$. In particular, $Z$ is
reduced. Notice also that $AB>0$, since the curves of $|L'|$ are $2$--connected by step
$2$, and thus $A$ and $B$ have nonempty intersection. Let $x_0\in A\cap B$, let
$y_0=\phi(x_0)$, and let $C$ be the pull--back of a general line through $y_0$; then
$C(2A+B)=L(2A+B)=4$ and  $x_0$ accounts at least for $3$ intersections. Thus
$\phi\inv(y_0)$ either consists of $x_0$ only, or contains also a point $x'_0$ that is not
a branch point of $\phi$; in either case $x_0$ is not a simple ramification point of
$\phi$ and therefore $Z$ is not smooth at $x_0$. Thus there is another line $R_1\subset Z$
that contains $x_0$. Write $\phi^*R_1= C_1=2A_1+B_1$. From $1=AL=2AA_1+AB_1$ we see that
either $AA_1=\frac{1}{2}$ and $AB_1=0$ or $AA_1=0$ and $AB_1=1$. On the other
hand, $A_1$  contains $x_0$ and thus we have $A_1B>0$ and $A_1A>0$. Thus we have a
contradiction, and $Z$ is not a union of lines. \par
We can now consider an irreducible component $Z'$ of $Z$ that is not a line and a general
tangent line $R$ to $Z'$. The curve $C_0=\phi^*(R)$ is reduced, but singular at some
point $x$. It moves in a base point free continuous system on $W$. Set
$h^{-1}(x)=\{x_1,x_2\}$ and let $C_0'=h^*C_0$. Notice that the map
$h\colon C_0'\to C_0$ is \'etale. Moreover $C_0'$ is singular at $x_1$ and $x_2$, and we
can  apply theorem (3.2) from \cite {ML}. Then we have $C'_0=A'+B'$ with
$A'$, $B'$ reduced and with no common component, since
$C'_0$ is reduced as well as $C_0$, and $A'B'=2$. Actually $A'\cap B' = \{x_1,x_2\}$,
which proves that $A'$ and $B'$ are smooth at $x_1$ and $x_2$.\medskip

\noindent {\it Step 4: One has $A'^2=B'^2=2$ hence $2A'$ and $2B'$ are numerically
equivalent to $L'$. }
Since $A'$ and $B'$ move without fixed components on $V$, we have $A'^2\geq 0$ and
$B'^2\geq 0$. Furthermore we have $A'^2+B'^2=4$ and $L'A'=A'^2+2$ and $L'B'=B'^2+2$.
Suppose $A'^2=0$, hence
$L'A'=A'B'=2$. We claim that in this case $A'$ is irreducible: in fact, if $A'=A_1+A_2$,
then $A_1^2$, $A_2^2\ge 0$, since $A_1$ and $A_2$ move, and thus $A_1^2=A_2^2=0$ and
$A_1\sim_{num}A_2$, $A_1B'=A_2B'=1$, contradicting the fact that the curves of $|L'|$ are
$2$--connected. Thus the general curve $A'$ is irreducible and moves in an irrational
pencil ${\cal A'}$ on $V$. The involution
$\iota$ determined by $h\colon  V\to W$ fixes $C_0'$, hence it maps $A'$ to an irreducible
curve $A''$. If $A'=A''$, then  there exists $A\subset C_0$ on $W$ such that
$A'=h^*A$, $A^2=0$, $LA=1$. Thus $A$ is smooth rational and, since $A'$ is general, $W$ is
covered by rational curves,  contradicting  
$\kappa(W)\ge 0$.  So $A''$ is contained in $B'$. The curve $A''$ also moves in an
irrational pencil ${\cal A}''$, and $A'A''\geq 2$, since $A'$ and $A''$ both contain $x_1$
and $x_2$.  Write $B=A''+D$, $C_0'=A'+A''+D$; since $C_0'A'=L'A'=2$, we get
$A'A''=2$ and $DA'=0$. Since $D$ also moves on $V$ without fixed components, it consists
of curves of ${\cal A'}$, hence $D^2=0$. Since $L'$ is fixed by $\iota$, we have
$L'A''=L'A'=2$, $A''^2=A'^2=0$, and therefore: $2=A''L'=A''(A'+A''+D)=2+A''D$, $A''D=0$.
Thus $D$ and
$A''$ are numerically equivalent, but this contradicts $A'A''=2$.
\par

Suppose that $A'^2=1$. By proposition (0.18) of \cite{CCML} 
we deduce that $A'$ is smooth irreducible and $V$ is isomorphic to the
symmetric product of $A'$. The canonical maps of
symmetric products are well known. Thus, the fact that the canonical map of $V$
is not birational, since it factors through $h$,  tells us that either
$3\ge p_a(A')=q(V)=g-1$ or $A'$ is hyperelliptic of genus $p_a(A')\geq 4$. The
former case is impossible by corollary \ref{quattrocor}. The latter case is
also impossible because $|L'|$ cuts out on $A'$ a base point free $g_3^1$.
Hence we are left with the only possibility $A'^2\geq 2$ and, similarly,
$B'^2\geq 2$, which implies the assertion.\par

\noindent {\it Step 5: the divisors $A'$ and $B'$ are exchanged by $\iota$.}
The divisor $\iota(A')=A''$ is contained in $C'_0$ and is numerically
equivalent to $A'$, since $2A'$ and $2A''$ are both numerically equivalent to
$L'$. If $A'=A''$, then there exists $A$ on $W$ such that $h^*A=A'$, $A^2=1$.
We apply proposition (0.18) of \cite{CCML}  to the pull-back of $A$ to the
minimal desingularization $\tilde{W}$ of $W$ and deduce  that $\tilde{W}$ is
birational to the symmetric product of $A$, contradicting $q(W)=0$. If  $A'$
is irreducible, this is enough to prove that  $A''=B'$. So assume that  $A'$
is  reducible and write $A'=N+M$, with $N, M$ effective nonzero. Then
$2=A'^2=A'N+A'M$, hence $A'N=A'M=1$ since $A$, as well as $L'$, is ample. This
proves that $N,M$ are both irreducible. Since they move on $V$, we have
$N^2\geq 0$ and $M^2\geq 0$ and the index theorem yields $N^2=M^2=0$, $NM=1$
and $N$ and $M$ both describe base point free pencils on $V$.  Since $A'\ne
A''$, $B'$ and $A''$ have at least a common component. Thus we may write
$B'=M'+N'$, where $M'$ is equal to, say, $\iota(M)$. 
We have $M'B'=MA=1$, ${M'}^2=M^2=0$, hence $B'N'=1$ and $N'$ is irreducible by the
ampleness of $B'$. If $\iota(N)=N'$, then the claim is proven. So assume $\iota(N)\ne
N'$. Then we have $\iota(N)=N$ and there exists $N_0\subset W$ such that $h^*N_0=N$. It
follows that $LN_0=1$,  and thus $N_0$ is a rational curve.
This is simpossible, since otherwise $W$ would be covered by rational curves.
 Thus $\iota(N)=N'$ and $\iota$
exchanges $A'$ and $B'$.\par

\noindent{\em Step 6: conclusion of the proof.}
We use the  notation introduced in step $1$.
By step 5, if the base point of the Albanese map $\alpha\colon V\to Prym(C',C)$ is
invariant for $i$, then $\alpha_*B'=(-1)_*\alpha_*A'$, since $\alpha$ is equivariant with
respect to $\iota$ on $V$ and multiplication by $-1$ on $Prym(C',C)$. Thus 
$\alpha_*(C'_0)=\alpha_*(A')+\alpha_*(B')$ and $2\alpha_*A'$ represent the same
cohomology class. By Welters criterion, this implies that
$\alpha_* A'$ is equivalent in cohomology to $\frac{1}{(g-2)!}\wedge^{g-2}\Xi$. By the
criterion of Matsusaka--Ran, $Prym(C'C)$ is isomorphic as a principally polarized abelian
variety either to a Jacobian or to a product of Jacobians. This contradicts step $1$, and
the proof is complete. \hfill\qed 

\begin{rem}\label{trix} {\em The same ideas we exploited in the proof of the
previous proposition would also yield the following result: {\em in case (ii)
of proposition \ref{classificagood}, the generating pair is either the pair of
the example  \ref {Rita} or it is obtained from Beauville's example
\ref{Beauville}, with a blow-up procedure}. We will next prove the same theorem
with a different technique, which also seems illuminating to us. Hence we give
here only an idea of its proof with the present methods.\par

Assume for simplicity $L^2=3$. Then one considers the finite, degree $3$ map $\phi\colon 
W\to {\pp^2}$ determined by $|L|$. First one shows that no line is in the branch divisor
of $\phi$. Then one proves the existence of a $1$-dimensional family of reduced curves
$C_0'\in |L|$ which split as $C_0=A+B$, with $AB=2$. This implies that $A^2=B^2=1$. At
this point one uses proposition (0.18) from \cite {CCML} and proves that $V$ is birational
to the symmetric product of $A=B$. The fact that the canonical map of $V$ is not
birational tells us that either $g\leq 3$, which leads to the two cases which actually
occur, or $A$ is hyperelliptic of genus $g\geq 4$. But this not possible because, via $h$,
$A$ is birational to the image $C_0$ of $C'_0$, and $C_0$ has a base point free $g_3^1$.} 
\end {rem}

The rest of this section is devoted to the analysis of case (ii) of
proposition \ref{classificagood} under the hypothesis  that $L$ be
ample. We prove the following result:
\begin{thm}\label{trigonalcase}Let $(h\colon V\to W, L)$ be a good generating
pair of genus $g$ such that $L$ is ample and $h^0(W,L)=3$.
Then $L^2=3$ and:\par
(i) either there exists a smooth plane quartic $\Gamma$ such that $(h\colon
V\to W, L)$ is constructed from $\Gamma$ as explained in example \ref{Rita};
\par
(ii) or $(h\colon V\to W, L)$ is obtained from Beauville's example
\ref{Beauville} via a simple blow-up of weight $1$.\end{thm} 
By propositions \ref{classificagood} and \ref{basepoint}, a pair
satisfying the assumption of theorem \ref{trigonalcase} either has $L^2=3$ and
$|L|$ is base point free or has $L^2=4$ and $|L|$ has a simple base point. So,
up to a simple blow up of the pair, we may assume that $L^2=3$ and $|L|$ is
base point free. Thus for the rest of the section we make the following
assumption:
\begin{assu}\label{trig}
$(h\colon V\to W, L)$ is a good generating pair of genus $g$ such that $L$ is ample,
$h^0(W,L)=3$,  $L^2=3$ and $|L|$ is base point free.
\end{assu}
If assumption \ref{trig} holds, then $|L|$ defines a finite morphism $f\colon
W\to\pp^2$ of degree $3$. The restriction of $f$ to  the general curve $C\in
|L|$ exhibits $C$ as a triple cover of $\pp^1$ showing that $C$ is trigonal. 
Given a curve $C$ of genus $g$, a degree $3$ map $f\colon C\to\pp^1$, and an
unramified double  cover $\pi\colon  C'\to C$, the {\it trigonal construction}
(\cite{Re}, cf. \cite{Beauspec}) yields a degree $4$ map $\phi\colon
D\to\pp^1$, where $D$ is a smooth curve of genus $g-1$ and $\phi$ has no
double fibre, such that the Jacobian of $D$ and $Prym(C',C)$ are isomorphic as
principally polarized abelian varieties. We briefly recall the trigonal
construction. One considers the induced morphism $\pi^{(3)}\colon  C'^{(3)}\to
C^{(3)}$  between the symmetric products of $C'$ and $C$. The curve
$\tilde{D}=\pi^{(3)^{-1}}(g_3^1)$ has a natural morphism $\tilde{D}\to \pp^1$;
it turns out  that $\tilde{D}\to\pp^1$  splits as the disjoint union of  two
isomorphic smooth connected  degree $4$ covers $\phi_i\colon D_i\to \pp^1$,
$i=1,2$, and one can set $D=D_1$, $\phi=\phi_1$.\par
The trigonal construction is a one--to--one correspondence, whose inverse is
the Recillas' construction (\cite{Re}, cf. \cite{LB} page 391). Given  a
smooth genus $g-1$ curve $D$ with a degree $4$ morphism
$\varphi\colon D\to\pp^1$ without double fibres, one defines a curve
$C'\subset  D^{(2)}$ by setting: $$C'=\{p_1+p_2\in D^{(2)}|  
p_1+p_2+p_3+p_4\,\, \hbox{is a fibre of $\varphi$  for some}\, \,p_1,  p_3\in
D\}.$$  The curve $C'$ is smooth and connected of genus $2g-1$, and has  a
natural  free involution $\sigma$, which maps an element
$p_1+p_2$ (in a fibre of $\varphi$) to the complementary element $p_3+p_4$. If
$\pi\colon C'\to C=C'/<\sigma>$\, denotes  the natural projection, it is easy to check that
$C$ is trigonal.\par

Recillas' correspondence has been generalized in \cite{cas2}, where the author
introduces the {\em discriminant } of a degree $4$ Gorenstein cover
$\varphi\colon Z\to Y$, which is a degree $3$ morphism $f\colon 
\Delta(Z)\to Y$. We recall that a cover
$\varphi\colon Z\to Y$ is said to be
a Gorenstein cover if the scheme theoretic fibre $\varphi^{-1}(y)$  is Gorenstein over
$k(y)$ for every $y\in Y$ (cf. \cite{cas1}). If $Y=\pp^2$, then the  discriminant
construction  gives a one--to--one correspondence between the following
objects:\par

(A) normal Gorenstein covers  $f\colon W\to \pp^2$ of degree $3$ such that the
singularities of
$W$ are at most RDP's and such that  there exists a double cover $h\colon V\to W$ branched
exactly over the singularities of $W$;\par

(B) degree $4$ Gorenstein covers $\varphi\colon Z\to \pp^2$ with $Z$  smooth
such that:
\begin{list} {(\roman{arab})}{\usecounter{arab}}
\item for every
$y\in \pp^2$ the Zariski tangent space to the fibre
$\varphi\inv(y)$ has dimension $\le 1$ at each point.
\item the set
$R_0\subset \pp^2$ of points $y$ such that  the fibre $\varphi\inv(y)$ is isomorphic
either to $\spec \, \C[t]/(t^4)$ or to
$\spec \, \C[t,s]/(t^2+1, s^2)$ is finite.
\end{list}
The properties of this correspondence ensure that the branch loci
of the associated covers $\varphi$ and  $f=\Delta(\varphi)$ coincide as
divisor of $\pp^2$. Moreover, the singularities of $W$ occur  precisely over 
the points $y\in R_0$. Notice that, in the case we are interested in, the
singular locus of $W$ is not empty (see the proof of theorem
\ref{prym}), and therefore $R_0$ is not empty. Finally, fibrewise, $Z_y$ is the
base locus of a pencil of conics whose discriminant is $W_y$.\par
Assumption \ref{trig} allows us to apply the trigonal contruction to the
present case. Thus, given a good generating pair $(h\colon V\to W,L)$ as in
\ref{trig},  there exists a unique degree
$4$ Gorenstein cover $\varphi\colon  Z\to \pp^2$ as in (B) such  that the morphism
$f\colon W\to
\pp^2$ associated to the system $|L|$ is obtained from $\varphi$ via the trigonal
construction. We denote by $|M|$ the pull-back to $Z$ of the linear system of  lines in
$\pp^2$.

\begin{lem}\label{isotrivial} The smooth elements of $|M|$ are isomorphic curves  of genus
$g-1\ge 2$. 
\end{lem} \proof Let $H$ be a general line in $\pp^2$, let
$D=\varphi\inv H$, let $C=f\inv H$ and let $C'\to C$ be the unramified cover determined by
$h$. By theorem \ref{prym},
the Prym variety $P=Prym(C',C)$ is independent of $H$. On the other hand,
$C'\to C$ is obtained from $D$ via the trigonal construction, and thus $P$ and the
Jacobian of $D$ are isomorphic as p.p.a.v.'s. In particular, since the genus of $|L|$ is 
at least $3$, the genus of
$|M|$ is at least $2$. By the global Torelli theorem for curves, the isomorphism class of
$D$ is also independent of
$H$. This implies that the natural map from the open set of smooth curves of
$|M|$ to the moduli space of curves is constant. \hfill\qed

\begin{lem}\label{nonsingular} Let $y\in \pp^2$ and let  $|M_y|\subset |M|$ be the 
pull-back on $Z$ of the pencil of lines through $y$. Then the general curve of $|M_y|$ is
smooth.
\end{lem} \proof The base scheme of $|M_y|$ is $\varphi\inv(y)$. The statement  follows by
Bertini's theorem since $\varphi\colon Z\to \pp^2$ satisfies condition (i) of (B).
\hfill\qed

\begin{lem}\label{ruled}
$Z$ is a minimal geometrically ruled surface,  and the smooth elements of
$|M|$ are sections of the ruling.
\end{lem} \proof Denote by $R\subset Z$ and by $B\subset \pp^2$ the
ramification divisor and the branch divisor of $\varphi$. By condition (B),
the ramification order of $\varphi$ along each component of $R$ is $\le 3$,
each component of $R$ is mapped birationally onto its image and different
components of $R$ are mapped to different components of $B$. Let $(M_t)_{t\in
\pp^1}$ be  a general  pencil contained in $|M|$ and assume that $M_0$ is
singular. By applying stable reduction, one can replace $M_0$ by a stable
curve $M'_0$.  Lemma \ref{isotrivial} implies that $M'_0$ is isomorphic to
$M_t$, for $t$ general. 

Assume that there exists a component $\Theta$ of $R$ such that
$\Delta=\varphi(\Theta)$ is not a line. If $\varphi$ is ramified of order $3$
along $\Theta$, then the inverse image $M_0$ of a generic line tangent to
$\Delta$  has an ordinary cusp over the tangency point and it is smooth
elsewhere. It follows that $M_0$ is irreducible. Since $p_a(M)>1$, $M_0$ is
not rational; the special fibre of the stable reduction of a general pencil
containing $M_0$ is the union of the normalization $M''_0$ of $M_0$ and of a
smooth elliptic curve meeting $M'_0$ at one point, but this is impossible by
the remark above. If $\varphi$ is simply ramified along $\Theta$, take 
$(M_t)$ to be the pull-back of a pencil of lines such that $M_0$ is the
pull-back of a line  simply tangent to $\Delta$ at a  point $y_0$  and 
meeting $B$ transversely elsewhere. Then $M_0$ has an ordinary node at a point
$x_0$ such that  $\varphi(x_0)=y_0$ and no other singularities. By the remarks
above, the curve $M_0$ is not  semistable; therefore we have  $M_0=M'_0+F$,
where $F$ is a smooth  rational curve, $M'_0$ is isomorphic to $M_t$ for $t$
general,  and $M'_0F=1$. We have: $4=M^2= M_t(M'_0+F)$, $M_tM'_0\geq 3$ (since
$M'_0$ is not hyperelliptic),  and thus  $M_tM'_0= 3$ and $F^2=0$. Noticing
that $y_0$ is a general point of $\Delta$, it follows that $Z$ is ruled. 
Since the system $|M|$ is ample, $M F=1$ and, by lemma \ref{isotrivial}, the
curves of $|M|$ are not rational, $Z$ is geometrically ruled and minimal.

So we have proven that either $Z$ is as claimed or all the  components of $B$ are lines.
Let $\Delta\subset B$ be a line; by condition (B), it is not possible that
$f^*\Delta=2A$. Thus $\varphi^*\Delta=mA+B$, with $m\le 3$, $A$ irreducible, $B$
nonempty and $A$ not contained in $B$. Then one can argue as in step $3$ of
the proof of proposition \ref{quattro} and prove that $B$ is not a union of
lines. \hfill\qed

\begin{prop}\label{quartic} Let $B$ be the base curve of the ruled surface
$Z$ of lemma \ref{ruled} and let $p\colon Z\to B$ be the projection. Then
there exists a birational morphism $s\colon B\to\Gamma \subset \pp^2$ such
that:\par

i) $\Gamma$ is either a smooth quartic or a quartic with a double point;\par

ii) $Z=\pp(s^*T_{\pp^2}(-1))$, and $p_*{\cal O}_Z(M)=s^*T_{\pp^2}(-1)$.
\end{prop} \proof According to lemma \ref{ruled}
there exists a rank $2$ bundle $E$ on $B$ such that $Z=\pp(E)$ and 
$p_*{\cal O}_Z(M)=E$ (in particular, $\deg(\det E)=4$). Let $D$ be a smooth
curve in $|M|$, which we may identify with $B$ via the map $p|_{D}$. Then
$M|_D$  is identified with $\det E$. By condition (B), (ii), if $D$ is
general, then $\varphi|_D$ has no multiple fibre, while if $\varphi(D)$
contains a point of $R_0$ (which, as we know, is not empty) then
$\varphi|_{D}$ has at least one multiple fibre.  So  the restriction of $|M|$
to $D$ is not a complete system, i.e. $h^0(B, \det E)=3$. Let $s\colon
B\to\pp^2$ be the morphism given by the linear system  $|\det E|$ and let
$\Gamma=s(B)$. If $\Gamma$ were a conic,  then  the map $\varphi|_D$ would
have two multiple  fibres  for every smooth $D$ of $|M|$, contradicting
condition (B). So $\Gamma$ is a quartic and $s$ is birational. Since $B$ has
genus $g-1\ge 2$, it follows that $\Gamma$ is either smooth or it has one
double point and  $\det E=s^*{\cal O}_\Gamma(K_\Gamma)$.\par

Let $U\subset H^0(Z, M)$ be the subspace such that 
$\pp(U)=\varphi^*|{\cal O}_{\pp^2}(1)|$.  If we identify $U$ with a subspace of
$H^0(B, E)$, then the natural sheaf map $U\otimes \OZ\to E$ is surjective
($|M|$ is base-point free). Moreover, the map $\wedge^2 U\to
H^0(B,s^*K_{\Gamma})$  is an isomorphism (this follows from the discussion
above, since we have shown that $|M|$ does not restrict to the same $g^1_4$ on
all the curves of $|M|$).

If we choose a basis for $U$, then we have a short exact sequence:
\begin{equation}\label{eulero} 0\to\OO_B(-s^*K_\Gamma)\to \OO_B^3\to E\to 0.
\end{equation}
Let the inclusion $\OO_B(-s^*K_\Gamma)\to \OO_B^3$ be given by $(s_0, s_1,
s_2)$, where $s_i\in H^0(\Gamma,K_\Gamma)$, 
$i=0,1,2$, and let $S$ be the subspace of
$H^0(\Gamma,K_\Gamma)$ spanned by $s_0,s_1,s_2$. Notice that 
$\dim S\ge 2$, since $E$ is torsion free. If $\dim S=2$, then it is clear that
$E=\OO_B\oplus \OO_B(s^*K_\Gamma)$  and condition ii) above is not satisfied. Thus $s_0,
s_1, s_2$ are independent and sequence (\ref{eulero}) is the pull-back via the map $s$ of
the twisted Euler sequence:

$$ 0\to\OO_{\pp^2}(-1)\to \OO_{\pp^2}^3\to T_{\pp^2}(-1)\to 0. \hfill\mbox{\qed}$$

\medskip Now we are ready to finish the proof of  theorem \ref{trigonalcase}:

\begin{prop} Notation as in proposition \ref{quartic}. The surface $Z$ is the
normalization of the incidence  surface
$Y=\{(p,l)\in \Gamma\times(\pp^2)^*|p\in l\}$, and the maps 
$p\colon Z\to B$ and 
$\varphi\colon Z\to\pp^2$ are induced by the projections of $Y$ onto
$\Gamma$ and $(\pp^2)^*$ respectively.\par

Let $f\colon W\to\pp^2$ be the triple cover obtained from $\varphi\colon Z\to\pp^2$ via
the  discriminant contruction, $h\colon V\to W$ the corresponding double cover and
$L=f^*{\cal O}_{\pp ^2}(1)$. Then:
\begin{list} {(\roman{rom})}{\usecounter{rom}}
\item if $\Gamma$ is smooth, then $V=Sym^2(\Gamma)$,  and $(h\colon
V\to W, L)$ is as in example \ref{Rita};
\item if $\Gamma$ has a double point, write $p+q$ for the only
effective divisor linearly equivalent to $s^*K_{\Gamma}\otimes K_B^*$.
Then:
\begin{list} {(\alph{alpha})}{\usecounter{alpha}}
\item
$V$ is the blow-up of $Sym^2(B)$ at $p+q$, namely it is the blow-up of the
Jacobian $J=J(B)=Pic^2(B)$ of $B$ at the points corresponding to $K_B$ and
$p+q$;
\item $W$ is obtained as the quotient of $V$ by the involution which is induced
on $V$ by the birational involution on $Sym^2(B)$ which associates to the
general divisor $x+y$ the divisor $|s^*K_\Gamma-x-y|$. Notice that $W$ is the
blow-up of the Kummer surface $Kum(J)$ at a smooth point;
\item the generating pair $(h\colon V\to W, L)$ is obtained from example
\ref{Beauville} by a simple blow-up of weight $1$.
\end{list}
\end{list}
\end{prop} \proof We keep the notation of the proof of proposition
\ref{quartic}.  Assume first $\Gamma$ is smooth. Then the first assertion
immediately follows by the well known fact that ${\pp}(T_{\pp^2}(-1))$ is the
incidence correspondence inside ${\pp^2}\times ({\pp^2})^*$. Having in mind 
Recillas' construction described at the beginning of this section, also part
i) immediately follows. The case $\Gamma$ singular is completely analogous and
can be dealt with in the same way. We leave the details to the reader.
\hfill\qed \section{The other cases}\label{altricasi}
In this section we collect some  information on  pairs that are not good or
not of Kodaira dimension $2$. We start by classifying  non good degree $2$
pairs with $L^2=4$. (We recall that by propositions \ref{kappa} and
\ref{Reider} such a pair always has $L^2\le 4$.)
\begin{prop}\label{classnongood} Let $(h\colon V\to W, L)$ be a non good 
generating pair of degree $2$ with $L^2=4$; then there exist smooth curves
$C_i$, $i=1,2$, of genus $g_i>0$ and double covers $\phi_i\colon C_i\to\pp^1$
such that $(h\colon V\to W, L)$ is obtained by a sequence of unessential 
blow--ups  from a generating pair constructed from $\phi_i\colon C_i\to\pp^1$
as in example \ref{product}. \end{prop} \proof  By proposition \ref{mintame}
we can assume that the pair is minimal.\par

\noindent Let $C\in |L|$ be general and let $C'=h^*C$. By \cite{Mumford}, p.
346, we see that the Galois group $G$ of the composition of $C'\to C$ with 
the hyperelliptic involution on $C$ can be identified with $\Z_2\times \Z_2$.
As in the proof of lemma \ref{bielliptic}, denote by $\sigma$ the element of 
$G$ such that $C'/<\sigma>=C$ and by $\sigma_i$ ($i=1,2$) the remaining non
trivial elements. For $i=1,2$, set $p_i\colon C'\to C_i=C'/<\sigma_i>$ the
corresponding projection and notice that $C_i$ is a smooth curve of genus
$g_i$, where $g_1+g_2=g-1$ and there exists a cartesian diagram:  
\begin{equation}\label{hyperell} \begin{array}{rcccl}
\phantom{1} &C'&\stackrel{\pi_2}{\rightarrow} &C_2&
\phantom{1} \\
\scriptstyle{\pi_1}\!\!\!\!\!\! & \downarrow & \phantom{1} &
\downarrow &
\!\!\!\!\!\!\scriptstyle{\phi_2}
\\
\phantom{1} & C_1 & \stackrel{\phi_1}{\rightarrow} &{\pp^1} &
\phantom{1}
\end{array}
\end{equation} where, for $i=1,2$, 
$\phi_i\colon C_i\to {\pp^1}$ is a double cover. In the present case there is
an isomorphism $Prym (C',C)=J(C_1)\times J(C_2)$ as  principally polarized abelian
varieties, and $A=Alb(V)$ is also isomorphic to $Prym (C',C)$ by theorem
\ref{prym}. We can assume that $g_1\le g_2$ and the condition that $C'$ is not
hyperelliptic ensures that
$g_1>0$. Notice the existence of commutative diagrams:\par
\begin{equation}\label{hyperelldia}
\begin{array}{rcccl}
\phantom{1} &C'&\stackrel{}{\rightarrow} &Prym(C',C)&
\phantom{1} \\
\scriptstyle{\pi_i}\!\!\!\!\!\! & \downarrow & \phantom{1} &
\downarrow &
\!\!\!\!\!\!\!\!\!\!\!\!\!\!\!\!\!\!\scriptstyle{p_i}
\\
\phantom{1} & C_i & \stackrel{}{\rightarrow} &{J(C_i)} &
\phantom{1}
\end{array}
\end{equation}
\noindent where $C'\to  Prym(C',C)$ is the Abel-Prym map, $C_i\to J(C_i)$ is
the Abel-Jacobi map  and
$p_i\colon  Prym(C',C)=J(C_1)\times J(C_2) \to J(C_i)$ is the $i$-th projection,
$i=1,2$.\par
As in the proof of lemma \ref{bielliptic}, one shows that there exist involutions
 $\tau_1$, $\tau_2$ on $V$ that act on   $C'$  as  $\sigma_1$, respectively, 
$\sigma_2$. Clearly, the involution $\iota$ associated to $h$ is equal to $\tau_1\circ
\tau_2$. We denote by
$S_i$ the quotient surface $V/\!\!<\tau_i>$, by $h_i\colon V\to S_i$ the projection onto
the quotient and by $C_i$ the image in $S_i$ of a general $C'$.  The
singularities of $S_1$ and $S_2$, if any, are $A_1$ points and $q(S_i)=g_i$.
More precisely, we claim that  $J(C_i)$ is the Albanese variety of $S_i$.
Indeed, the map $V\to J(C_i)$ obtained by composing the Albanese map of $V$
with the projection onto $C_i$ is equivariant with respect to $\tau_i$,
provided that the base point of the Albanese map is invariant for $\tau_i$.
Thus we have an induced map $S_i\to J(C_i)$ and thus the Albanese variety
$A_i$ of $S_i$ is isogenous to $J(C_i)$. To show that this isogeny is actually
an isomorphism, it is enough to remark that the map $H_1(V, \Z)\to
H_1(J(C_i),\Z)$ is surjective, since it is the composition  of   $H_1(V,\Z)\to
H_1(A,\Z)$, that is an isomorphism up to torsion, and of $H_1(A,\Z)\to
H_1(J(C_i),\Z)$ which is surjective. On the other hand,  $H_1(V, \Z)\to
H_1(J(C_i),\Z)$ is also the composition of $H_1(V, \Z)\to H_1(S_i,\Z)$ and
$H_1(S_i, \Z)\to H_1(J(C_i),\Z)$, hence the latter map is surjective and $A_i$
is isomorphic to $J(C_i)$.\par
We claim that $S_i$ is birational to $\pp^1\times C_i$. Indeed,  by
proposition \ref{xiaofib} the curve $C_i$ does not vary in moduli  and the
Albanese image of $S_i$ is a curve. By proposition \ref{xiaose}, this
concludes the proof of the claim. In particular the Albanese image of $S_i$ is
the curve $C_i$.\par 
Composing $h_i$ with the Albanese map $S_i\to C_i$, we  get morphisms
$f_i\colon V\to C_i$, $i=1,2$. Denote by $F_i$ a fibre of $f_i$. The Index
theorem applied to $F_1+F_2$ and $L'$ gives $2(F_1F_2)L'^2\le
[L'(F_1+F_2)]^2=16$, namely $L'^2\le 8$, $L^2\le 4$. If $L^2=4$, then 
$F_1F_2=1$ and $L'$ is numerically equivalent to $2F_1+2F_2$. Thus
$f=f_1\times f_2\colon  V\to C_1\times C_2$ is birational, and therefore it is
an isomorphism since $V$ is minimal. One has: $\tau_1=\sigma_1\times Id$,
$\tau_2=Id \times \sigma_2$, $\iota=\sigma_1\times\sigma_2$ and the curves of
$|L'|$ are  invariant for $\tau_1$, $\tau_2$ and it is easy to see that
$(h\colon V\to W, L)$ is precisely as in example \ref{product}. \hfill\qed

\medskip
Next we classify  pairs of degree $2$ and Kodaira dimension $0$. 
\begin{prop}\label{kappa0} Let $(h\colon  V\to W,L)$ be a  generating pair of degree
$2$ and genus $g$;  if the Kodaira dimension of the pair is $0$, then it can be obtained
from example \ref{Beauville} (Beauville's example) by a sequence of simple blow-ups, only
three of which at most essential, of weight $1$.
\end{prop} \proof Assume that the pair has Kodaira dimension $0$ and is minimal. By
proposition \ref{regular} and corollary \ref{pg>0}, we see that $g=3$ and the irregularity
of $V$ is $2$, hence $V$ is  an abelian surface. Since $q(W)=0$ by proposition
\ref{regular},
$W$ is the Kummer surface of $V$. By theorem \ref{prym}, if $C\in |L|$ is general then $V$ is
isomorphic to $Prym(C',C)$ and thus, in particular, it is principally polarized. In
addition,
by Welters criterion, $C'$ is a divisor of type $(2,2)$ and thus we 
have precisely example \ref{Beauville}. By corollary \ref{mintame}, this implies that if the
pair is not minimal, then it is obtained from example \ref{Beauville} by a sequence of
blow--ups of weight $0$ or $1$. Since $L$ is big by assumption and  example
\ref{Beauville}  has $L^2=4$, there are at
most $3$ blow--ups of weight $1$ in the sequence. \hfill\qed
\medskip
The next result is an almost complete  classification of  pairs of degree $2$ and 
Kodaira dimension $1$.

\begin{prop} \label{kappa1} Let  $(h\colon V\to W,L)$ be a  generating
pair of degree $2$ and genus $g$ with Kodaira dimension $1$. Then  there exist
an elliptic curve $E$ and an hyperelliptic curve $B$ of genus $g-2\ge 2$ such
that $(h\colon V\to W,L)$ is obtained by a sequence of simple blow--ups of
degree $0$ or $1$ from  one of the following:\par
(a) the pair constructed from $E$ and $B$ as in example \ref{product}. In this
case the pair is not good; \par
(b) a pair $(h_0\colon V_0\to W_0, L_0)$ such that $g=4$ (and thus $B$ has
genus $2$), $V_0=B\times E$ and $h_0\colon V_0\to W_0$ is the quotient map for
the $\Z_2$--action given by $(b,e)\mapsto (j(b),\sigma(e))$, where $j$ is the
hyperelliptic involution of $B$ and $\sigma$ is an involution of $E$ with
rational quotient. In this case  $L^2=2$, and, if the pair is good, then
$h^0(W,L)=2$. \end{prop} \proof  By  corollary \ref{mintame} we may assume
that $(h\colon V\to W,L)$ is minimal.\par
\noindent Let $V\to B$ be the elliptic fibration. By a result of Beauville
(see \cite {Beau4}, pg. 345) and by corollary \ref{albsurface}, $V$ is a
product   $B\times E$, where $E$ is the general fibre of $V\to B$. The
involution $\iota$ determined by $h$ on $V$ preserves the fibration $V\to E$,
and, since the quotient of $V$ by $\iota$ is regular, it acts on $B$ as an
involution $j$ with rational quotient. Thus $\iota$ can be written as
$(b,e)\mapsto (j(b), \sigma_b(e))$, where $\sigma_b\colon E\to E$ is an
automorphism of $E$. If $\sigma_b$ is a translation for every $b\in B$, then
the pull-back on $V$ of the nonzero $1$-form of $E$ is invariant for $\iota$,
but this contradicts the regularity of $W$. So $\sigma_b$ is not a
translation. Since $\iota^2=1$, one has $\sigma_b\circ\sigma_{j(b)}=1$, and,
if $b_0$ is a fixed point of $j$, then  $\sigma_{b_0}^2=1$.  So $\sigma_{b_0}$ acts on
$H^0(E,\omega_E)$ as multiplication by $-1$.  Since the possible actions of an
automorphism of $E$ on
$H^0(E,\omega_E)$ are a finite number, it follows that $\sigma_b$ acts on
$H^0(E,\omega_E)$ as multiplication by $-1$ for every $b\in B$. Thus we have 
$\sigma_b^2=1$, namely  $\sigma_b=\sigma_{j(b)}$. So $b\mapsto\sigma_b$ descends to a
map $B/<j>=\pp^1\to Aut(E)$ and it is therefore constant. Notice that $B$ has genus
$g-2\ge 2$, since $\kappa(V)=1$.\par Denote by $F$ the general fibre of the
pencil $p_1\colon W\to {\pp^1}=B/<j>$;
$F$ is isomorphic to $E$ and $p_1$ has $2g-2$ double fibres, each containing
$4$ nodes of $W$. Now, with the usual notation, we take $C\in |L|$ a general
curve and $C'=h^*C$.  Note that $CF$ is even, since $p_1$ has double fibres
and  the general $C$ contains no singular point of $W$. So we set $C F=2l$.
The system $|K_W|$ is equal to $|(g-3)F|$, and thus the adjunction formula on
$V$ gives: $$4(g-1)=C'^2+C'K_V= C'^2+4l(g-3).$$
If $l=1$, then we have $C'^2=8$, namely $L^2=4$, and $p_1$ restricts to a
$g^1_2$ on $C$, so that $C$ is hyperelliptic and the pair is not good. Thus 
proposition \ref{classnongood} implies that we are in case (a). Then we have
$0<C'^2= 4[g(1-l)+3l-1]$, which leaves us with the only possibility $l=2$,
$g=4$, $C'^2=4$, and therefore $C^2=2$  and this is case (b). If the pair is
good, then $h^0(W, L)=2$  by proposition \ref{classificagood}. \hfill\qed

\medskip

By  proposition \ref{regular}, only non good generating pairs can have
degree $3$ or $4$ and only for very restricted values of the genus $g$. The following
theorem gives some more information on this case: 

\begin{prop}\label{d=3} Let $(h\colon V\to W, L)$ be a  generating
pair of degree $d>2$ and genus $g$.  Then one of the
following holds:\par
i)  $d=3$, $q(V)=g=2$, $\kappa(V)\ge 0$,  $p_g(V)\ge 1$;\par
ii) $d=3$, $q(V)=4$, $g=3$,  $p_g(V)\geq 4$,  $V$ is a surface of general type
and its Albanese image is a surface;\par
iii) $d=4$, $q(V)=3$, $g=2$, $p_g(V)\geq 2$ and the Albanese image of $V$ is a
surface.  \end{prop} \proof The possible values of $d$ and $g$ and the
corresponding values of the irregularity of $V$ are given in proposition
\ref{regular}, as well as the assertion on the dimension of the Albanese image
of $V$ in case iii). The claim on the dimension of the Albanese image of $V$
in the case ii) follows by lemma \ref{isolemma}. So the Albanese image of $V$
is a surface, and thus $V$ is not ruled in these cases.\par 

Assume now that we are in case ii), so that $\kappa(V)\ge 1$, by the Kodaira--Enriques
classification of surfaces. If $\kappa(V)=1$, then the minimal
model $\bar{V}$ of $V$ is equal to $E\times B$, where 
$E$ is an elliptic curve and $B$ is a smooth curve of genus $3$, since otherwise the
Albanese image of $V$ would be a curve (cf. \cite{deb}, Lemma page 345), contradicting
theorem \ref{prym}. \par
We recall that $W$ is regular by proposition \ref{regular}
and that the surfaces $V$ and $W$ have the same Kodaira dimension by
proposition \ref{kappa}. Let
$\bar{p}\colon W\to \pp^1$ be the elliptic fibration on $W$. Clearly we have a commutative
diagram: $$\begin{array}{rcccl}
\phantom{1} &V&\stackrel{p}{\rightarrow} &B&\phantom{1} \\
\scriptstyle{h}\!\!\!\!\!\! & \downarrow & \phantom{1} &\downarrow &
\!\!\!\!\!\scriptstyle{\bar{h}}\\
\phantom{1} & W & \stackrel{\bar{p}}{\rightarrow} & \pp^1 &\phantom{1}
\end{array}$$
where $p\colon V\to B$ is the composition of $V\to\bar{V}$ and of the projection
$\bar{V}=B\times E\to B$, and 
$\bar{h}\colon B\to \pp^1$ has  degree $3$. The  map
$h\colon V\to W$, being finite, is obtained from $\bar{h}$ by base change and
normalization. Let
$y\in\pp^1$ be such that $\bar{h}$ is branched at $y$ and assume that $\bar{
h}^*y=2x_1+x_2$, with $x_1\ne x_2$. Since $h$ is \'etale in codimension $1$, in
particular it is not ramified   along $p\inv x_1$. It follows that
the fibre of $\bar{p}$ over
$y$ must be a double fibre. But then 
the fibre of $p$ over $x_2$ is a double fibre, since the diagram is commutative and
$\bar{h}$ is unramified at $x_2$.  This is impossible, since  $p$ is obtained from the
projection $B\times E\to B$  by a composition of blow--ups. Thus the ramification points of
$\bar{p}$ all have order $3$, and there are $5$ of them by the Hurwitz formula. By the
classical Riemann construction, the covering $\bar{p}\colon B\to\pp^1$ is determined, up
to isomorphism, by the branch points
$y_1,\ldots y_5\in \pp^1$ and  by permutations $\sigma_i\in S_3$ describing the local
monodromy at $y_i$. The $\sigma_i$ satisfy $\sigma_1\ldots\sigma_5=1$. Up to
renumbering the $y_i$, the only possibility is that there is a $3$-cycle $\sigma\in
S_3$ such that
$\sigma_i=\sigma$ for  $i=1\ldots 4$ and $\sigma_5=\sigma\inv$. By \cite{ritaabel},
proposition 2.1, there exists a cyclic cover with these properties, and thus
$\bar{h}\colon B\to\pp^1$ is cyclic, and the same is true for $h\colon V\to W$. 
The Galois group $\Z_3$ of $h$ acts also on the minimal model $\bar{V}=E\times B$ of
$V$,  and the quotient is a surface $\bar{W}$ with rational singularities.  The
minimal resolution $\tilde{W}$ of $\bar{W}$ has invariants $p_g=3$, $q=0$.  Arguing as
in the proof of proposition
\ref{kappa1}, one shows that a generator
$\gamma$ of $\Z_3$ acts on $V$ by $(b, e)\mapsto (\gamma b,\sigma_b e)$, where the
action of
$\gamma$ on $B$ is the one associated to the Galois cover $\bar{h}\colon V\to\pp^1$ and
$\sigma_b$ is an automorphism of order $3$ of $E$ that is not a translation. The
action of $\sigma_b$ on $H^0(E, \omega_E)$ is independent of $b$.  
Each of the curves $\{x_i\}\times E$, $i=1\ldots 5$,  contains $3$ fixed points of the
$\Z_3$ action on $\bar{V}$, and these are the only fixed points.
The surface $\tilde{W}$ has an $A_2$ singularity at the image of a fixed point $P$ of
$V$ if the representation of $\Z_3$ on the tangent space at $\bar{V}$ in $P$ is
contained in $SL(2,\C)$ and has a singulairty of type 
$\frac{1}{3}(1,1)$ otherwise. From the above description of the $\Z_3$ action, it follows
that $\bar{W}$ has either
$12$ points of type $A_2$ and $3$ points of type $\frac{1}{3}(1,1)$, or it has $3$
points of type $A_2$ and $12$ points of type $\frac{1}{3}(1,1)$. The Euler characteristics
of $\bar{V}$ and of $\tilde{W}$ are related by the formula:

$$\chi(\bar{V})=3\chi(\tilde{W})-\frac{1}{3}\alpha-\frac{2}{3}\beta$$
where $\alpha$ is the number of singularities of type $\frac{1}{3}(1,1)$ and $\beta$ is
the number of singularities of type $A_2$. Thus we have either $\chi(\tilde{W})=3$ or
$\chi(\tilde{W})=2$, contradicting  $\chi(\tilde{W})=4$. So this case does not occur, and we
have shown that
$\kappa(V)=2$ if $q(V)=4$. In particular, one has $p_g(V)\ge 4$.\par

We turn now to the case $q(V)=g=2$. Suppose that $V$ is ruled and denote by $p\colon V\to
B$ the Albanese pencil of $V$, where $B$ is a  curve of genus $2$. We can consider a
minimal model $\bar{V}$ of $V$ that is a geometrically ruled surface over $B$.  Arguing
exactly as in the previous case, one shows  that the Galois group of $h$ is isomorphic
to $\Z_3$. Therefore $\Z_3$ acts also on $\bar{V}$ with $8$ fixed points, giving rise to
$4$ singularities of type $A_2$ and $4$ singularities of type $\frac{1}{3}(1,1)$ of
the quotient surface $\bar{W}$. The minimal desingularization of $\bar{W}$ is a rational
surface and thus we have a
contradition, using again the formula above for the Euler characteristics. In
particular, $\chi(V)\ge 0$, namely
$p_g(V)\ge q(V)-1=1$. \hfill\qed

\bigskip

\noindent Authors' adddresses:\bigskip

\noindent Ciro Ciliberto and Francesca Tovena\par

\noindent Universit\`a di Roma Tor Vergata \par

\noindent Dipartimento di Matematica\par

\noindent Via della Ricerca Scientifica\par

\noindent 00133 Roma, Italia\bigskip

\noindent Rita Pardini\par

\noindent Universit\`a di Pisa \par

\noindent Dipartimento di Matematica "L. Tonelli"\par

\noindent Via F. Buonarroti 2\par

\noindent 56127 Pisa, Italia\medskip
\end{document}